\documentclass[]{arXivclass}

\usepackage{endnotes}
\let\footnote=\endnote
\let\enotesize=\normalsize
\def\notesname{Endnotes}%
\def\makeenmark{\hbox to1.275em{\theenmark.\enskip\hss}}
\def\enoteformat{\rightskip0pt\leftskip0pt\parindent=1.275em
  \leavevmode\llap{\makeenmark}}
\usepackage{natbib}
 \bibpunct[, ]{(}{)}{,}{a}{}{,}%
 \def\bibfont{\small}%

\TheoremsNumberedThrough     %
\ECRepeatTheorems

\EquationsNumberedThrough    %

\usepackage[hypertexnames=false]{hyperref}

\usepackage{mathtools, nccmath}
\usepackage{braket} %
\usepackage{siunitx}
\usepackage{enumitem}
\setlist[enumerate]{align=left}
\usepackage{amsmath,amssymb}
\usepackage{mathabx}
\usepackage{cleveref}
\usepackage{xcolor}
\usepackage{verbatim}
\usepackage[short]{optidef}
\DeclareDocumentEnvironment{miniinf}{D||{\defaultProblemFormat} O{\defaultConstraintFormat} D<>{} m m m m}{%
\ifthenelse{\equal{#3}{b}}{%
    \ifthenelse{\equal{#1}{s}}%
    {\setFormatShort{inf}{#4}\BaseMiniStar{#2}{#4}{#5}{#7}{inf}{#3}}%
    {\setFormatLong{infimum}{#4}\BaseMiniStar{#2}{#4}{#5}{#7}{infimum}{#3}}%
}{%
    \ifthenelse{\equal{#1}{s}}%
    {\setFormatShort{inf}{#4}\BaseMini{#2}{#4}{#5}{#6}{#7}{inf}}%
    {\setFormatLong{infimum}{#4}\BaseMini{#2}{#4}{#5}{#6}{#7}{infimum}}%
}%
}%
{\endBaseMini\toggletrue{bodyCon}}
                        
\usepackage{booktabs}
\usepackage{graphicx}
\usepackage{subfig} %
\usepackage{multirow}
\usepackage[tworuled,linesnumbered]{algorithm2e}

\newcommand{\gap}{\vspace{0.05in}}

\newcommand{\epc}{\hspace{1pc}}

\newcommand{\onebld}{{\bf 1}}

\newcommand{\wh}{\widehat}

\begin{document}

\TITLE{Classification and Treatment Learning with \\ 
Constraints via Composite Heaviside Optimization: \\
a Progressive MIP Method\footnote{Original date: \today} }

 \ARTICLEAUTHORS{%
\AUTHOR{Yue Fang}
\AFF{School of Management and Economics,  The Chinese University of Hong Kong, Shenzhen, China, \EMAIL{\tt fangyue@cuhk.edu.cn}} %
\AUTHOR{Junyi Liu}
\AFF{Department of Industrial Engineering, Tsinghua University, Beijing 100084, China, \EMAIL{\tt junyiliu@tsinghua.edu.cn}}
\AUTHOR{Jong-Shi Pang}
\AFF{Department of Industrial and Systems Engineering, University of Southern California, Los Angeles \EMAIL{\tt jongship@usc.edu}}
} %
\RUNTITLE{Classification and treatment learning  
by Heaviside composite optimization}

\RUNAUTHOR{Fang, Liu, and Pang}  
\vspace{4in}

\ABSTRACT{This paper proposes a Heaviside composite 
optimization approach and presents a progressive (mixed) 
integer programming (PIP) method for solving multi-class 
classification and multi-action treatment problems with 
constraints.  A Heaviside composite function is a composite 
of a Heaviside function (i.e., the indicator function of either 
the open $( \, 0,\infty )$ or closed 
$[ \, 0,\infty \, )$ interval) with a 
possibly nondifferentiable function.   Modeling-wise, we show 
how Heaviside composite optimization provides a unified formulation 
for learning the optimal multi-class classification and multi-action 
treatment rules, 
subject to rule-dependent constraints stipulating a variety of 
domain restrictions.   A Heaviside composite function has an 
equivalent discrete formulation 
and the resulting optimization problem can in principle be 
solved by integer programming (IP) methods.  Nevertheless, for
constrained learning problems with large data sets, 
a straightforward application of off-the-shelf IP solvers 
is usually ineffective in achieving global optimality.  To alleviate such 
a computational burden, our major contribution is the proposal of 
the PIP method by leveraging the effectiveness 
of state-of-the-art IP solvers for problems of modest sizes.
We provide the theoretical advantage of the PIP method 
with the connection to continuous optimization and show that the 
computed solution is locally optimal for a broad class of Heaviside 
composite optimization problems. 
The 
numerical performance of the PIP method is demonstrated by 
extensive computational experimentation.}

\SUBJECTCLASS{Integer programming, nonlinear programming.}


\maketitle

\section{Introduction}

\vspace{0.5ex}

Classification is a well-accepted statistical tool for solving
numerous applied problems in science, engineering, 
medicine, health, and economics disciplines.  As early as 
the mid 1960's (see \cite{Mangasarian65}), it is known
that linear programming can provide a very powerful computational 
platform for binary classification problems.  With leading work by
Breiman and his collaborators \cite{BFOStone84} that is exemplified 
by their classification and regression trees (CART), decision 
tree methods are widely used for classification problems.  
While a key strength of the tree methods is its  
well-interpretation, traditional methods for growing decision 
trees are mostly heuristic and myopic without optimality 
guarantees, and most restrictively, lack the ability to deal 
with problems with domain knowledge that leads to constraints
on the classification process.  In this unconstrained framework,
in order to learn optimal trees more rigorously and with 
guarantees, dynamic programming and mixed integer programming
(MIP) have recently received growing attention.  The latter 
approach is of particular interest due to the availability 
of off-the-shelf solvers for broad classes of MIP problems. 
\cite{BertsimasDunn17} propose the first 
MIP formulation based on the big-M technique for learning 
optimal classification trees with numerical results 
showing that the state-of-the-art Gurobi solver is applicable 
for such models with better performance than CART.  For 
classification tree problems with integer or categorical 
features, the more recent article (see \cite{AghaeiGomezVayanos21}) 
proposes a core integer programming formulation without the big-M 
technique, showing that it has a tighter LP relaxation and is 
amenable to Bender's decomposition. As reported in the last
reference, for problems with depths up 
to 4 and thousands of data points, Gurobi normally takes at 
least three to four hours to solve the corresponding MIP to 
optimality, while it could find high-quality solutions in
shorter time, though no optimality properties can be 
guaranteed if terminated early. 


\gap

Formulated in the continuum space of covariates, typically 
a finite-dimensional Euclidean space and thus distinct from 
the space of bivariate covariates that is the setting of the work by 
\cite{AghaeiGomezVayanos21}, covariate-dependent rule-based 
predictive treatment learning has been of increasing interest in 
several areas, including 
medical treatments, policy learning, and machine learning.  
In these applications, the treatment rules are used to define 
performance quantities to be optimized or domain-desirable targets 
to be satisfied.  Examples include population 
welfare (see \cite{athey21}) and the Gini index of income 
inequality in \cite{Gini1909} for 
policy learning (see \cite{YueFang23}), treatment outcome in 
health care (see \cite{QiCuiLiuPang19}), accuracy at the top 
(see \cite{BoydCortesMohriRadov12}) and its variations 
in \cite{AdamMachaSmidl23}  as well as goals such as fairness, 
recalls, and regrets in machine learning 
(see \cite{CHGWNYSridharan19,CotterJiangSridharan19,
ESMGSElidan17}).   The goal of optimization is to 
judiciously select the parameters of some prescribed 
families of classification (or score) function to optimize 
an overall performance objective (not only the count of
misclassifications) subject to domain requirements under target 
quantities, such as those targets just mentioned.  Indeed, the ability to
handle constraints is a distinguished advantage of an 
optimization approach over a heuristic approach. Moreover, by
considering and optimizing the physical or other welfare outcomes 
of the treatment rule, such as a patient's wellness after a
particular kind of treatment, we arrive at a class of 
multi-action treatment
problems that broaden a traditional classification problem
in which the outcomes concern primarily misclassification.

\gap

Due to the discrete nature of classification and treatment rules, 
the rule-dictated objective and constraint
functions are necessarily discontinuous; in fact, they can
often be represented as compositions of univariate indicator 
functions of 
1-dimensional intervals (the well-known Heaviside functions)
with the action-inspired classification functions, which may be 
nonsmooth but nevertheless continuous functions.  Such composite 
functions (with constant muplipliers) generally take the form: 
$\displaystyle{
\sum_{j=1}^J
} \, \psi_j \, \onebld_{[ \, 0,\infty \, )}( \phi_j(x) )$, 
for some positive integer $J$, constants $\{ \psi_j \}_{j=1}^J$,
and functions $\phi_j : \mathbb{R}^n \to \mathbb{R}$, where 
$\onebld_{[ \, 0,\infty \, )} (t) \triangleq \begin{cases} 
0 & \mbox{ if } t< 0  \\ 1 & \mbox{ if } t \geq 0  
\end{cases}$ is the indicator function of the closed intervale
$[ \, 0, \infty )$.  We call such composite functions
\emph{Heaviside composite functions} and an optimization problem 
defined by them a \emph{Heaviside composite optimization problem}
(HSCOP) which is formally stated as follows:
\begin{equation} \label{eq:HSCOP}
\displaystyle{
\operatornamewithlimits{\mbox{\bf maximize}}_{x \, \in \, 
X_{\rm AHC}}
} \ \theta(x) + f(x), \epc \mbox{where } \ 
f(x) \, \triangleq \, \displaystyle{
\sum_{j=1}^{J_0}
} \, \psi_{0j} \, \onebld_{[ \, 0,\infty \, )}( \phi_{0j}(x) ),
\end{equation}
where $\theta : \mathbb{R}^n \to \mathbb{R}$ is a Bouligand 
differentiable function \cite[Chapter~4]{CuiPang2021}
and 
\begin{equation} \label{eq:AHC set}
X_{\rm AHC} \, \triangleq \, \left\{ \, x \, \in \, P \, 
\left| \, \displaystyle{
\sum_{j=1}^{J_i}
} \, \psi_{ij} \, \onebld_{[ \, 0,\infty \, )}( \phi_{ij}(x) ) 
\, \geq \, b_i, \ i \, = \, 1, \cdots, m \, \right. \right\}.
\end{equation}
for a given polyhedral subset $P$ of $\mathbb{R}^n$, some 
positive integers $m$ and $\{ J_i \}_{i=0}^m$, constants 
$\{ b_i \}_{i=1}^m$,
nonnegative constants $\left\{ \{ \psi_{ij} \}_{j=1}^{J_i} 
\right\}_{i=0}^m$, and functions 
$\phi_{ij} : \mathbb{R}^n \to \mathbb{R}$. Providing a unified 
mathematical formulation and being the computational workhorse
for both treatment and classification 
problems with constraints, the topic of Heaviside composite
optimization has been the subject of the recent 
references such as \cite{CuiLiuPang22-piecewise,HanCuiPang23}, where 
the modeling breadth of the Heaviside functions in 
optimization is detailed.  The approach in these two studies
follows the tradition of continuous optimization where 
concepts of stationarity (i.e., necessary conditions for 
local optimality) are defined and surrogation methods 
\cite[Chapter~7]{CuiPang2021} for nonconvex, 
nondifferentiable methods provide the basis for computing such 
stationary solutions.   
\gap 

From the discrete perspective, an HSCOP 
can be reformulated as a mixed-integer program by employing binary
variables. Realizing that the straightforward application of 
off-the-shelf IP solvers 
is usually ineffective for constrained learning problems with large 
data sets, the present work takes a different route in addressing 
the HSCOP for classification and treatment learning.  Specifically, 
leveraging advances in integer programming but not relying on it in 
a straightforward way, we propose a methodology 
that bridges domains of continuous and integer 
programming. There are several distinct benefits of 
the proposed method, which we term the \emph{progressive
integer programming (PIP) method}.  To begin, starting at a 
feasible solution of the HSCOP, the method generates a sequence 
of feasible solutions with non-deteriorating objective values.
The computed solution at the termination of the method can be shown
to be a local optimizer for a broad  class of HSCOPs,
including those arising from the classification 
and treatment learning applications.  On the contrary, the 
latter locally optimizing property cannot be claimed by an IP 
solver applied to the HSCOP with the full set of integer 
variables when it is terminated prematurely before 
a globally optimal solution is found.  In practice, 
high-quality suboptimal solutions
(with guarantees of bounds from global optimality in 
verifiable cases) can be computed in a reasonable amount of 
time for problems that present challenges for the commercial 
IP solver Gurobi to even find a feasible solution.

\gap

The organization of the rest of the paper is as follows.  
In the next section, we introduce a multi-action treatment 
learning problem with the Gini fairness constraint from observational 
data, and the parameterized policy class for multiple actions without 
ambiguity by using margins. We provide the HSCOP formulation for 
treatment learning problems, highlighting the semicontinuity of the 
formulations which plays a key role for rigorous analysis and the 
algorithmic development to be presented subsequently 
({\bf Contribution~I}).  
In Section~\ref{sec:standard classification}, we discuss several
multi-class classification paradigms: standard versus the 
Neyman-Pearson (NP) paradigms with score-based and tree-based 
classification rules, with the NP framework naturally inducing
a constraint.
We show that these frameworks can also be formulated as HSCOPs
({\bf Contribution~II}).  From the perspective of 
continuous optimization, we summarize some known local 
stationarity results for an HSCOP in Section~\ref{sec:HSCOP}, 
which prepares for the presentation of the progressive integer 
programming (PIP) method subsequently in
Section~\ref{sec:the PIP method}  ({\bf Contribution~III}).  
Indeed, the local results are responsible for the important 
termination property of the PIP method mentioned above.
Finally, in the last Section~\ref{sec:numerics}, we present results
of an extensive set of computational experiments on the
treatment learning problem with the Gini constraint 
demonstrating the practical performance of the PIP method and 
its superiority over the Gurobi solver for problems with a large 
number of Heaviside terms ({\bf Contribution~IV}).

\gap

As a side note, we mention that the key idea
underlying the PIP method can be applied to many nonconvex
(non-Heaviside) 
optimization problems.  Some preliminary results conducted
outside of the present work have demonstrated
that the method is able to improve a given stationary solution
of a nonconvex quadratic program
and can compute high-quality solutions with properties
similar to those of the problems presented herein.  Additionally,
many other nonconvex optimization problems can similarly be solved. 
We believe that this is a major 
computational breakthrough in the area of nonconvex optimization outside 
the realm of global optimization.  
Details of all the extended developments are presently in progress.

\section{Multi-action Treatment with Income Fairness Constraint}
\label{sec:Treatment with Income}

We build up the core HSCOP model 
(\ref{eq:treatment optimization}) for treatment learning of this 
section by first 
describing the model framework of multi-action treatment with some key 
assumptions for accounting counterfactuals with observational data.  
We introduce the key model elements for optimization which include the 
income fairness constraint utilizing Gini measure under 
the \emph{Inverse probability weighting} (IPW) estimations 
(\cite{hirano2003efficient, kitagawa2018should, kitagawa2021equality}) 
and the parameterized policy class for assigning the treatment from 
multiple choices of action.  With all the background preparations, 
we then formally define the multi-action treatment learning optimization 
problem under the income fairness constraint
and present its workhorse formulation as an HSCOP
for whose solution the PIP method is designed.

\subsection{The probabilistic framework of multi-action treatment} 

There are $J$ treatments to be assigned according to an 
individual's covariate, which is a random vector 
$X : \Omega \to \mathbb{R}^p$ that is defined 
on the probability triple $( \Omega,{\cal F},{\cal P} )$, 
with $\Omega$ being the sample space, ${\cal F}$ the 
$\sigma$-algebra defined on $\Omega$, and
${\cal P}$ the probability operator defined on ${\cal F}$.  Let
$D: \Omega \rightarrow [ J ] \triangleq \{1, \ldots, J\}$
denote a random treatment assignment; we denote the realizations of 
$X$ and $D$ by $\xi \in \mathcal{X}$ and $d \in [ J ]$, respectively,
where $\mathcal{X}$ is a subset of $\mathbb{R}^p$.
Associated with each treatment $j \in [ J ]$ is a random outcome 
$Y(j) : \Omega \to \mathbb{R}$ whose support we assume is a bounded 
interval $[ \, 0, M \, ]$ for a given scalar $M > 0$.  Of interest
is the composition $Y \circ D$, written $Y(D)$ that is the random 
outcome under
the treatment assignment $D$.  The significance of this composite
random variable is that for an individual with the covariate $X$, 
one can observe the realization of $Y(D)$ from
experiments but not that of $Y(j)$ for $j$ outside the
observed range of the treatment $D$.  The latter outcomes
are referred as \emph{counterfactuals}.  A \emph{treatment policy} 
is a mapping $g: \mathcal{X} \to [J]$ which maps an 
individual's covariate vector $\xi \in \mathcal{X}$ to a treatment 
$j \in [ J ]$.  In addition to the bounded outcome assumption,
the following are two standard assumptions in the
causal inference literature:

\gap

\noindent {\bf (A1)} Conditional exchangeability (Unconfoundedness):
for every treatment $j \in [ J ]$, the outcome $Y(j)$ is 
independent of the treatment $D$ given the covariate $X$; 
written as $Y(j) \perp D \mid X$. 

\gap

\noindent {\bf (A2)}] Strict overlap: there exists  
$\kappa\in ( \, 0, 1/2 \, )$, such that for every treatment 
$j \in [ J ]$, the {propensity score}
$e_j(\xi) \triangleq \mathbb{P}\left[ D=j\mid X = \xi \right] \in 
[ \, \kappa, 1-\kappa \, ]$ for all $\xi \in \mathcal{X}$.

\gap

\noindent The former assumption means 
that after conditioning on realized characteristics, the potential 
outcomes are independent of the treatment assignments. Formally, 
by Assumption (A1), the expected welfare and the cumulative 
distribution function (cdf) of the outcome under policy $g$ 
are given, respectively, by
\[
\begin{array}{cll}
&&\mathbb E\left[ \, \mathbb E\left[ \, Y(g(X)) \, | \,
X \, \right] \, \right]  =  \displaystyle{
\sum_{j=1}^J  
} \, \mathbb E\left[ \, \onebld\left\{ g(X) = j \right\} \, 
\mathbb E\left[ \, Y(j) \, | \, X \, \right] \, \right] 
  =   \displaystyle{
\sum_{j=1}^J
} \, \mathbb{E}\left[ \, \onebld\left\{ g(X)=j \right\} \,
\mathbb{E}\left[ \, Y(D) \mid X; D=j \, \right] \, \right] \\ [0.2in]
&&F_g(y) = \displaystyle{
\sum_{j=1}^J
} \, \displaystyle{
\int_{\, \mathcal X}
} \, \left\{ \, \onebld\left\{ g(\xi)=j \right\} \, 
F_{Y(j)\mid \xi}(y) \, \right\} 
\, d \, \mathbb{P}_{X}(\xi)
\end{array}
\]
where $\onebld\{ \bullet \}$ is the true-false indicator function,
which is equal to one if the statement is true and zero otherwise. 
Based on these two fundamental quantities, the rest of the
derivation of the core treatment model proceeds in the several
steps.

\gap

\noindent $\bullet $ {\bf Inverse probability 
weighting (IPW) estimator:}
The joint probability distribution of the triple
$\left\{ X,D,Y(D) \right\}$ is generally unknown.  In the 
empirical approach to the treatment problem of interest, 
we suppose that independently and
identically distributed (iid) data triples 
$\{ \, ( X^s, D^{\, s}, Y_s ) \, \}_{s=1}^N$ are available 
from observational experiments; in particular, each sample $X^s$ is
a copy of the random variable $X$ that corresponds to treatment
assignment $D^{\, s}$ with outcome $Y_s = Y(D^{\, s})$.  A difficulty 
in treatment learning is the estimation of 
the conditional distribution $F_{Y(j)|\xi}$ 
and the conditional expectation 
$\mathbb E\left[ \, Y(j) \, | \, X=\xi \, \right]$ 
from only the observable samples of the treatment-outcome pair 
$( D,Y(D) )$ given the observed covariate $\xi$. 
Based on such samples,
IPW provides a way to estimate the counterfactual outcomes 
for all subjects in the population that were assigned the given
treatments.  

For each covariate $\xi \in \mathcal{X}$, let
${\cal S}_{\xi} \triangleq \{ s \in [ N ] : X^s = \xi \}$,
and thus $\displaystyle{
\bigcup_{\xi \in \mathcal{X}}
} \, {\cal S}_{\xi} = [ N ]$.
Note that each sample $X^s$ is associated with a unique 
$\xi \in \mathcal{X}$, but each $\xi \in \mathcal{X}$ 
may correspond to multiple samples, reflecting that the data may 
contain the outcomes of multiple treatments applied to the same 
covariate.  Furthermore, for each $s \in N$, we denote 
$j_s \triangleq D^{s}$ as a unique integer in $[ \, J \, ]$.  
Note, however, that for two different samples 
$s \neq s^{\, \prime}$, we may have 
$(X^s,j_s) = (X^{s^{\prime}},j_{s^{\prime}})$, yet
$Y_s \neq Y_{s^{\prime}}$, reflecting that the outcomes of
the same treatment applied to the same covariate may be different
in the data, due perhaps to external factors. Based on the samples 
$\{(X^s,D^s,Y_s)\}_{s=1}^N$, let $\wh{e}_j(\xi)$ be an estimation 
of the propensity score $e_j(\xi)$. 
By reformulating the expected welfare 
and the cdf of welfare under the policy $g$ as follows,
\[
\begin{array}{rll}
\mathbb{E}_g[ \, Y \, ] & = & \displaystyle{
\sum_{j=1}^J 
}\, \mathbb E\left[ \, \frac{Y(j) \, \onebld\{D=j\}}{e_j(X)}
\, \onebld\left\{ g(X) = j \right\} \, \right]\end{array} \]
\[
\begin{array}{rl} 
F_g(y) & =  1 - \mathbb E_{X}\Big[\, \displaystyle{
\sum_{j=1}^J
}  \onebld\{ g(X) = j \} \mathbb{P}\left( Y(j) > y \, | 
X \right) \Big]
\\[0.1in] 
&  = 1 - \mathbb E_{X}\Big[ \, \displaystyle{
\sum_{j=1}^J
} \, \displaystyle{
\frac{1}{e_j(X)}}
\mathbb{P}\left( Y(D) > y, D=j \mid X \right) \, \Big] \\[0.1in]
&=  1 - \mathbb E_X\Big[ \, \displaystyle{
\sum_{j=1}^J
} \, \displaystyle{
\frac{\mathbf{1}(D=j)}{e_j(X)}
} \, \onebld(Y>y) \, \onebld\left\{ g(X)=j \right\} 
\, \Big],
\end{array} \]
we can obtain the IPW-based estimator of expectation and cdf 
of the outcome variable $Y$ under the policy $g$ as follows,
\[ \begin{array}{ll}
\widehat{\mathbb E}^{\, \rm IPW}_g [ \, Y \, ] \, & = \,
\displaystyle{
\frac{1}{N}
}  \displaystyle{
\sum_{s=1}^N
} \, \displaystyle{
\frac{Y_s}{\wh{e}_{j_s}(X^s)}
} \, \onebld\{ g(X^s) = j_s \}   =
\displaystyle{
\frac{1}{N}
}  \displaystyle{
\sum_{\xi \in \mathcal{X}}
} \ \displaystyle{
\sum_{s \in {\cal S}_{\xi}}
} \, \displaystyle{
\frac{Y_s}{\wh{e}_{j_s}(x)}
}  \onebld\left\{ g(\xi) = j_s \right\},\\[0.2in]
\wh{F}^{\, \rm IPW}_g(y) & = 1- \displaystyle{
\frac{1}{N}
} \, \displaystyle{
\sum_{s=1}^N
} \, \displaystyle{  
\frac{\onebld ( Y_s > y  )}{\wh{e}_{j_s}(X^s)}
} \, \onebld\left\{ g(X^s) = j_s \right\} =1- \displaystyle{
\frac{1}{N}
} {\displaystyle{
\sum_{\xi \in \mathcal{X}}
} \displaystyle{
\sum_{s \in {\cal S}_{\xi}}
}}  \displaystyle{  
\frac{\onebld ( Y_s > y  )}{\wh{e}_{j_s}(\xi)}
} \, \onebld\left\{ g(\xi) = j_s \right\}.
\end{array} \]

Doubly robust (DR) estimation (see \cite{Scharfstein99}) 
is another well-known approach for  accounting causal effects, which 
combines IPW with estimated propensity score and expected outcome for 
each treatment arm given covariates so that the overall estimate 
is consistent if either one is. Though the DR approach is mostly 
applied for estimating expected outcome, it is also applicable to
the Gini coefficient (to be defined momentarily) through a little 
more involved derivations.  Moreover, the treating learning problem 
under the DR estimation shares the same Heaviside composite optimization 
structure as the one under IPW estimation, except with different 
coefficients.  Thus the succeeding
HSCOP formulations and the proposed PIP algorithm are all applicable 
to both estimation approaches.  For the sake of brevity, we omit the 
details of the DR estimator in the rest of the paper.

\noindent $\bullet $ {\bf The Gini constraint:}
For a treatment learning problem that is implemented for social 
welfare, it is important that policies maintain a satisfactory 
level of social inequality within social groups.  We adopt the 
utilitarian utility as objective function, represented by the mean 
income of the population of interest, and employ the Gini coefficient 
as the measure of income inequality. Various methods exist for 
measuring income inequality, as discussed in economics 
(see \cite{Schutz1951,deMaio2007}); among such measures, 
Gini coefficient is one of the 
most widely used.  It is derived from the Lorenz curve, which
illustrates the percentage of total income earned 
by cumulative percentage of the population.  In a scenario of perfect 
equality, the Lorenz curve coincides with the 45 degree line, 
symbolizing the equality line. The Lorenz 
curve diverges from this line as inequality increases. 
The Gini coefficient is computed by 
dividing the size of the area between Lorenz curve and the equality 
line, by the size of the area under the equality line.  A higher 
Gini coefficient corresponds to a greater degree of inequality.  
Although mean-log deviation, coefficient of variation, and decile 
ratios are alternative measures of 
inequality, we opt for the Gini coefficient due to its widespread
recognition.  Moreover, due to the generality of HSCOP, the methods 
in this present paper are able to accommodate many other types 
of inequality measures as well.

\gap

In terms of the IPW-based cdf 
$\wh{F}^{\, \rm IPW}_g$ and expected outcome 
$\widehat{\mathbb E}^{\, \rm IPW}_g [ \, Y \,]$, given 
a prescribed threshold $\alpha > 0$, the empirical Gini coefficient 
upper bounded by $\alpha$ is expressed by the inequality
\begin{equation} \label{eq:Gini constraint}
1 - \displaystyle{
\frac{\int_{\, 0}^{\, M}(1 - \wh{F}^{\, \rm IPW}_g(y))^2dy}{
\wh{\mathbb{E}}^{\, \rm IPW}_g(Y)}
} \, \leq \, \alpha,
\end{equation}
which is equivalent to
\[
\begin{array}{l}
\displaystyle{
\int_{\, 0}^{\, M}
} \, \wh F_g^{\, \rm IPW}(y)^2dy + 
(1+\alpha) \, \wh{\mathbb{E}}_g^{\, \rm IPW}(Y) - M \, \geq \, 0
\\ [0.2in]
\ \Longleftrightarrow \ \displaystyle{
\frac{1}{N^2}
} \, \displaystyle{
\sum_{s=1}^N
} \, \displaystyle{
\sum_{t=1}^N
} \, \displaystyle{
\int_{\, 0}^{\, M}
} \, \left( \, \displaystyle{
\frac{\onebld\{g(X^s) = j_s\}}{\wh{e}_{j_s}(X^s)}
} \, \right) \, \left( \, 
\displaystyle{
\frac{\onebld\{g(X^t) = j_t\}}{\wh{e}_{j_t}(X^t)}
} \, \right) \, \onebld ( Y_s \leq y  ) \onebld ( Y_t \leq y  ) dy  
\\  [0.2in]
\hspace{0.4in} + \, \displaystyle{
\frac{1 + \alpha}{N}
} \, \displaystyle{ 
\sum_{s=1}^N
} \, \displaystyle{
\frac{Y_s}{\wh{e}_{j_s}(X^s)}
} \, \onebld\{g(X^s) = j_s\} \, \geq \, M  \\ [0.2in] 
\Longleftrightarrow \ \displaystyle{
\frac{1}{N^2}
} \, \displaystyle{
\sum_{s=1}^N
} \, \displaystyle{
\sum_{t=1}^N
} \, \displaystyle{
\frac{M - \max\{Y_s, Y_t\}}{\wh{e}_{j_s}(X^s) \, \wh{e}_{j_t}(X^t)} 
} \, \onebld\{g(X^s) = j_s\} \, \onebld\{g(X^t) = j_t\} \\ [0.2in] 
\hspace{0.4in} + \, \displaystyle{
\frac{1 + \alpha}{N}
} \, \displaystyle{
\sum_{s=1}^N
} \, \displaystyle{
\frac{Y_s}{\wh{e}_{j_s}(X^s)}
} \, \onebld\{g(X^s) = j_s\} \, \geq \, M .
\end{array}
\]
$\bullet $ {\bf The parameterized policy class:}  
For computational purposes, we restrict the class of 
treatments rules to affine mappings with parameters to be determined 
by optimization.  Specifically, let 
$\{ \xi \mapsto \xi^{\top} \beta^j + b_j \}_{j=1}^J$ be a family 
of linear functions, where each $\beta^j \in \mathbb{R}^p$ 
is an unknown class-dependent parameter
and $b_j \in \mathbb{R}$ is a known scalar considered as the
``base score'' of class $j$ (for example, $b_j$ could be the base
benefit of treatment $j$ such as the negative of its cost that is
associated with a welfare maximization problem). 
A basic policy rule is defined by $g(\xi) \in \displaystyle{
\operatornamewithlimits{\mbox{\bf argmax}}_{j \in [ J ]}
} \, ( \xi^{\top} \beta^{\, j} + b_j )$.  With this rule, we have
\begin{equation} \label{eq:policy rule}
\onebld\left\{ g(\xi) = j \right\} \, = \, 
\onebld_{[ \, 0,\infty )}(h_j(\xi,\boldsymbol{\beta})), \epc
\xi \, \in \mathcal{X},
\end{equation}
where with 
$\boldsymbol{\beta} \triangleq ( \, \beta^{\, j} )_{j=1}^J$,
\begin{equation} \label{eq:PA score}
h_j(\xi,\boldsymbol{\beta}) \, \triangleq \, 
\xi^{\top} \beta^{\, j} + b_j - \displaystyle{
\max_{m \neq j}
} \, ( \, \xi^{\top} \beta^{\, m} + b_m \, ), \epc j \, \in \, [ J ].
\end{equation}
One has to be mindful about the 
mathematical rigor of the identity in (\ref{eq:policy rule}) 
that would be in jeopardy when $\displaystyle{
\operatornamewithlimits{\mbox{\bf argmax}}_{j \in [ J ]}
} \, ( \xi^{\top} \beta^{\, j} + b_j )$ is not a singleton, unless
we allow the policy map $g$ to be multi-valued.  Although such
non-uniqueness happens rarely in practice, especially in the presence
of numerical errors in computation, there are multiple 
ways to resolve the resulting ambiguity in defining the policy $g$.
One is by natural selection of the one with the largest welfare
(assuming that there is no tie with the welfare among the treatments).
Another way is by random assignment among the maximizing labels with
equal welfare.
A third way is to modify the score function $h_j$ with
a scalar $\tau_j > 0$; namely, via the alternative classification
function:
\begin{equation} \label{eq:h hat}
h_j^{\tau}(\xi,\boldsymbol{\beta}) \, \triangleq \, 
\xi^{\top} \beta^{\, j} + b_j - \displaystyle{
\max_{m \neq j}
} \, ( \, \xi^{\top} \beta^{\, m} + b_m \, ) - \tau_j,
\end{equation}
which leads to the policy:
\begin{equation} \label{eq:tau policy}
\onebld\left\{ g^{\tau}(\xi) = j \right\} \, = \, 
\onebld_{[ \, 0,\infty )}( \, h_j^{\tau}(\xi,\boldsymbol{\beta}) 
\, ), \epc \xi \, \in \mathcal{X},
\end{equation}
that can be contrasted with (\ref{eq:policy rule}).
Clearly, for $\tau_j > 0$, there exists at most 
one $j$ such that $h_j^{\tau}(\xi,\boldsymbol{\beta}) \geq 0$,
thus ensuring the well-definedness of the 
policy (\ref{eq:tau policy}).
While this modification resolves the uniqueness of the treatment 
class $j$, it raises the 
possibility of non-existence of a treatment $j\in[J]$ such that 
$h_j^{\tau}(\xi,\boldsymbol{\beta}) \geq 0$ for a given parameter 
$\boldsymbol{\beta}$ and a realization $\xi$ of the randomness.  
To guard for such non-existence, we introduce 
a \emph{least-residual} formulation of the empirical Gini constraint 
in the next subsection; see (\ref{eq:treatment optimization}).
Note that the function $h_j^{\tau}$ is just a constant shift of the 
function $h_j$ in (\ref{eq:PA score}), thus preserves the same 
functional properties.  

\subsection{The treatment optimization problem as a HSCOP}

Putting together the above modeling elements, we introduce the 
following IPW-based optimization formulation of the treatment 
assignment problem of maximizing welfare  
subject to the Gini constraint (\ref{eq:Gini constraint}): for given
weights $\lambda \geq 0$ and $\rho > 0$, 
\begin{equation} \label{eq:treatment optimization}
\begin{array}{l}
\displaystyle{
\operatornamewithlimits{\mbox{\bf maximize}}_{\gamma \, \geq 0; 
\, \boldsymbol{\beta}}
} \ \underbrace{\displaystyle{
\frac{1}{N}
} \, \displaystyle{
\sum_{\xi \in \mathcal{X}}
} \ \displaystyle{
\sum_{s \in {\cal S}_{\xi}}
} \, \displaystyle{
\frac{Y_s}{\wh{e}_{j_s}(\xi)}
} \, \onebld_{[ \, 0,\infty )}(h_{j_s}(\xi,\boldsymbol{\beta}))}_{
\mbox{empirical welfare}}
- \lambda \, \underbrace{\displaystyle{
\sum_{j=1}^J
} \, R_j(\beta^{\, j})}_{\mbox{regularizer}}
- \rho \, {\gamma}
\\ [0.5in]
\mbox{\bf subject to}  \quad 
\displaystyle{
\frac{1}{N^2}
} \, \displaystyle{
\sum_{\xi \in \mathcal{X}}
} \, \displaystyle{
\sum_{s \in {\cal S}_{\xi}}
} \, \displaystyle{
\sum_{\xi^{\, \prime} \in \mathcal{X}}
} \, \displaystyle{
\sum_{t \in {\cal S}_{\xi^{\, \prime}}}
} \, \underbrace{\displaystyle{
\frac{M - \max\{Y_s, Y_t\}}{\wh{e}_{j_s}(\xi) \, 
\wh{e}_{j_t}(\xi^{\, \prime})} 
}}_{\mbox{nonnegative}} \, \underbrace{
\onebld_{[ \, 0,\infty )}(h_{j_s}(\xi,\boldsymbol{\beta})) \, 
\onebld_{[ \, 0,\infty )}(h_{j_t}(\xi^{\, 
\prime},\boldsymbol{\beta}))}_{\mbox{
$= \, \onebld_{[ \, 0,\infty )}( \min( 
h_{j_s}(\xi,\boldsymbol{\beta}), \, 
h_{j_t}(\xi^{\, \prime},\boldsymbol{\beta})) )$}} \\ [0.45in]
\epc \qquad \qquad  + \, \displaystyle{
\frac{1 + \alpha}{N}
} \ \displaystyle{
\sum_{\xi \in \mathcal{X}}
} \, \displaystyle{
\sum_{s \in {\cal S}_{\xi}}
} \, \left( \, \displaystyle{
\frac{Y_s}{\wh{e}_{j_s}(\xi)}
} \, \onebld_{[ \, 0,\infty )}(h_j(\xi,\boldsymbol{\beta})) 
\, \right) +{\gamma} \, \geq \, M,
\end{array} \end{equation}
where we have employed (a) the regularizers $R_j(\beta^{\, j})$
to induce sparsity, with a prominent example being
the $\ell_1$-function: 
$R_j(\beta^{\, j}) = \| \, \beta^{\, j} \, \|_1$; and (b)
the scalar variable ${\gamma}$
as the residual of the Gini constraint to guard against its 
infeasibility, such that the least-residual Gini 
formulation with $\gamma = 0$ implies that the Gini constraint 
(\ref{eq:Gini constraint}) is satisfied (empirically) 
and $\gamma > 0$ for the opposite.  It is easy to see that
(\ref{eq:treatment optimization}) is an instance of the
HSCOP (\ref{eq:HSCOP}); the coefficients $\{ \psi_{ij} \}$ 
associated with the Heaviside composite functions are indeed
all nonnegative; moreover, the original definition
(\ref{eq:PA score}) or its $\tau$-modfication (\ref{eq:h hat})
renders each $h_j$ concave and piecewise affine. 
%
%
In general, with these properties, the objective function of 
(\ref{eq:HSCOP}) is upper semicontinuous
and the feasible set $X_{\rm AHC}$ is closed.  As documented in
the two recent references \cite{CuiLiuPang22-piecewise,HanCuiPang23},
these distinguished properties of the problem have an important role 
to play in its theoretical analysis and practical solution. 

\section{Related Classification Problems} 
\label{sec:standard classification}

Before discussing
further properties of the problem (\ref{eq:HSCOP}),
we present in this 
section several types of classification problems that can be 
formulated as a constrained HSCOP to illustrate its extensive applications.  The
resulting formulations will highlight that some care is required
to obtain desirable semicontinity properties to facilitate
the computational solution of the problems.

\subsection{Two paradigms of classification: standard and
Neyman-Perason}

Let $(X,Y)$ be a pair of jointly distributed random variables with 
$X$ denoting a $p$-dimensional attribute vector, and 
$Y \in [ \, J \, ]$ denoting the corresponding label from 
$J$ classes.  In standard classification problems, the probability
of correct classifications $\mathbb P(g(X)=i \mid Y=i)$ 
and the probability of misclassifications, 
$\mathbb P(g(X)=j \mid Y=i)$ for $j \neq i$ are two common 
criteria to be optimized, the former to be maximized and the 
latter minimized.   With $\mathcal G$ denoting a hypothesis class 
of mappings from attributes to labels, and based on the criterion
of maximizaing the probability of correct 
classification, we obtain the problem:
\begin{equation} \label{eq:standard classification}
\displaystyle{
\underset{g \in \mathcal G}{\mbox{\bf maximize}}
}  \sum_{i \in [J]} \mathbb{P}\left(g(X)=i \mid Y=i \right).
\end{equation}
In real-world applications such as medical 
diagnosis and cybersecurity, it is crucial to recognize that some 
types of classification errors are much more severe than others. 
The Neyman-Pearson (NP) paradigm is a common 
framework to incorporate asymmetric control over classification 
errors, but most of the current study is limited to binary 
classification problems where only type-1 error and type-2 error 
are considered.  To accommodate the 
asymmetric control over misclassification errors, we let 
$\mathcal E_1 \subseteq \{(i,j): i \neq j\}$ denote a set of 
misclassification label pairs that lead to more severe 
consequences and thus are controlled within the constraints 
by a prescribed threshold, and let 
$\mathcal E_2 \subseteq \{(i,j): i \neq j\}$ with 
$\mathcal E_1 \cap \mathcal E_2 = \emptyset$ and 
$\mathcal E_1 \cup \mathcal E_2 = \{(i,j): i \neq j\}$
denote the complementary set of misclassification 
label pairs for which the misclassification errors 
are to be minimized.  By letting ${\cal G}$ 
denote a hypothesis class of mappings from attributes to 
classified labels, we may formulate the population-level 
multiclass Neyman-Pearson classification learning problem 
as the following affine chance-constrained stochastic program
with a threshold $\gamma > 0$:
\begin{equation} \label{eq:NP classification}
\begin{array}{ll}
\displaystyle{
\operatornamewithlimits{\mbox{\bf minimize}}_{
g \in \mathcal G}
} & \displaystyle{
\sum_{(i,j) \in \mathcal E_2}
}  \, w_{ij} \, \mathbb P(g(X)=j \mid Y=i) \\ [0.25in]
\displaystyle{
\operatornamewithlimits{\mbox{\bf subject to}}
} & \displaystyle{
\sum_{(i',j') \in \mathcal E_1}
} \, w_{i^{\, \prime} j^{\, \prime}} \, \mathbb P(g(X)=
j^{\, \prime} \mid Y=i^{\, \prime}) \, \leq \, \gamma,
\end{array}
\end{equation} 
where the family $\{w_{ij}\}$ consists of (nonnegative) 
weights associated with the misclassification errors,
whose separation in the objective function and constraints 
reflects the different attitudes towards the errors 
within the two groups.  For instance, the weight $w_{ij}$ 
can be set as an estimation of the probability 
$\mathbb P(Y=i)$.  When 
$\mathcal E_1 = \emptyset$ and $w_{ij} = \mathbb P(Y=i)$, 
the above NP problem reduces to the standard (mis)classification 
problem. 

\gap

To obtain computationally tractable formulations of 
(\ref{eq:standard classification}) and
(\ref{eq:NP classification}), we employ a parameterized
family of classification functions 
$\{ h_j(\xi,\boldsymbol{\beta}) \}_{j=1}^J$ such that 
$\onebld\{g(X)=j\} = 
\onebld_{[\, 0, \infty)}(h_j(X, \boldsymbol{\beta}))$; moreover, 
for a given batch of data samples $\{(X^s, Y_s)\}_{s=1}^N$ of  
attribute-label pairs we approximate
\begin{equation} \label{eq:probability by Heaviside}
\mathbb{P}(g(X) = j \mid Y = i ) \, = \,  
\displaystyle{
\frac{\mathbb{E}[ \onebld(g(X)=j) \,
\onebld(Y = i ) ]}{\mathbb{P}(Y = i )}
} \, \overset{\rm SAA}{\approx \approx} \,
\displaystyle{
\frac{1}{| {\cal S}_i|}
} \, \displaystyle{
\sum_{s \in {\cal S}_i}
} \, \onebld_{[ \, 0,\infty )}(h_j(X^s,\boldsymbol{\beta})),
\end{equation}
where ${\cal S}_i = \{s: Y_s=i\}$.   Substituting this approximation
into (\ref{eq:standard classification}) and adding the appropriate 
regularizer $R(\boldsymbol{\beta})$, we obtain the following
unconstrained HSCOP formulation:
\begin{equation} \label{eq:HSCOP standard classification} 
\displaystyle{
\underset{\boldsymbol{\beta}}{\mbox{\bf maximize}}
} \ \displaystyle{
\sum_{i \in [J]} 
} \, \displaystyle{
\frac{1}{| {\cal S}_i|}
} \, \displaystyle{
\sum_{s \in {\cal S}_i}
} \, \onebld_{[ \, 0,\infty )}(h_i(X^s,\boldsymbol{\beta})) 
- \lambda R(\boldsymbol{\beta}).
\end{equation}
In contrast, if we 
substitute the approximation (\ref{eq:probability by Heaviside}) 
directly into (\ref{eq:NP classification}), we would obtain
a minimization problem lacking lower semicontinuity.  So by  
allowing margins in the classification and also replacing the 
closed Heaviside function with the open Heaviside function, we
obtain the following NP classification problem: 
into (\ref{eq:NP classification}), we obtain
\[ \begin{array}{ll}
\displaystyle{
\operatornamewithlimits{\mbox{\bf minimize}}_{\boldsymbol{\beta}}
} & \displaystyle{
\sum_{(i,j) \in \mathcal E_2}
}  \, \displaystyle{
\frac{w_{ij}}{| {\cal S}_i|}
} \, \displaystyle{
\sum_{s \in {\cal S}_i} 
} \, \onebld_{( \, 0,\infty )}\left( h_j(X^s,\boldsymbol{\beta}) + \tau_j 
\right) + \lambda R(\boldsymbol{\beta}) \\ [0.25in]
\displaystyle{
\operatornamewithlimits{\mbox{\bf subject to}}
} & \displaystyle{
\sum_{(i^{\, \prime},j^{\, \prime}) \in \mathcal E_1}
} \, \displaystyle{
\frac{w_{i^{\, \prime} j^{\, \prime}}}{| {\cal S}_{i^{\, \prime}} |} 
} \, \displaystyle{
\sum_{s \in {\cal S}_{i^{\, \prime}}} 
} \, \onebld_{( \, 0,\infty )}\left( 
h_{j^{\, \prime}}(X^s,\boldsymbol{\beta})+ \tau_{j^{\, \prime}}  \right) 
\, \leq \, \gamma,
\end{array}
\]
Since 
$\onebld_{( \, 0,\infty )}(t) = 1 - \onebld_{[ \, 0,\infty )}(-t)$,
the above minimization problem is equivalent to the maximization
problem below:
\[ \begin{array}{ll}
\displaystyle{
\operatornamewithlimits{\mbox{\bf maximize}}_{\boldsymbol{\beta}}
} & \displaystyle{
\sum_{(i,j) \in \mathcal E_2}
}  \, \displaystyle{
\frac{w_{ij}}{| {\cal S}_i|}
} \, \displaystyle{
\sum_{s \in {\cal S}_i} 
} \, \onebld_{[ \, 0,\infty )}\left( -h_j(X^s,\boldsymbol{\beta}) 
- \tau_j \right) + \lambda R(\boldsymbol{\beta})\\ [0.25in]
\displaystyle{
\operatornamewithlimits{\mbox{\bf subject to}}
} & \displaystyle{
\sum_{(i^{\, \prime},j^{\, \prime}) \in \mathcal E_1}
} \, \displaystyle{
\frac{w_{i^{\, \prime} j^{\, \prime}}}{| {\cal S}_{i^{\, \prime}} |} 
} \, \displaystyle{
\sum_{s \in {\cal S}_{i^{\, \prime}}} 
} \, \onebld_{[ \, 0,\infty )}\left( 
-h_{j^{\, \prime}}(X^s,\boldsymbol{\beta}) - \tau_{j'} \right) \, \geq \, 
\displaystyle{
\sum_{(i^{\, \prime},j^{\, \prime}) \in \mathcal E_1}
} \, \displaystyle{
\frac{w_{i^{\, \prime} j^{\, \prime}}}{| {\cal S}_{i^{\, \prime}} |} 
} - \gamma,
\end{array}
\]
which is of the same form as (\ref{eq:HSCOP}); the only noteworthy
difference is that when $h_j$ is given by
(\ref{eq:PA score}), $-h_j(X^s,\bullet)$ is convex while 
remaining piecewise affine.

\subsection{Tree-based classification rule}  
\label{subsec:tree-based classification}
\vspace{0.5ex}

In what follows, we show that optimal tree classification 
problem in the reference \cite{BertsimasDunn17} can also 
be formulated as a HSCOP with the same structure as the previous 
two problems.  Consider a standard classification tree problem 
with $J$ classes without asymmetric error control.  As before, 
let $\{ (X^s, Y_s) \}_{s=1}^N$ be the 
data set of size $N$; let $\mathcal T_\ell$ denote the set of 
leaf nodes, and $\mathcal T_{\mathcal B}$ denote the set 
of branching nodes.  The tree is constructed 
by splitting each branch node $k \in {\cal T}_{\cal B}$ into 
a right branch: $( a^k )^{\top}x \geq b_k$ and a left branch 
$( a^k )^{\top}x \leq b_k - \varepsilon$
for a given constant $\varepsilon > 0$ (as a replacement for 
strict inequality), where each
such pair $( a^k,b_k ) \in \mathbb{R}^{p+1}$ is to be 
determined. 
For each leaf node 
$t \in \mathcal T_\ell$, let $A_R(t)$ denote the set of 
ancestors of $t$ whose right branch has been followed on 
the path from the root node to $t$, and similarly for 
$A_L(t)$. Such a leaf node $t \in \mathcal T_\ell$ is assigned to a class $j_t \in [ J ]$. 
Under this setting, the tree classification problem is 
formulated as the following problem of maximizing the 
total number of correctly classified data points over the 
selection of the parameter tuple 
$( \mathbf{a},\mathbf{b} ) \triangleq 
\{(a^k, b_k)\}_{k \in {\cal T}_{\cal B}}$ and the 
leaf node assigning rule $\{ j_t \}_{t \in {\cal T}_{\ell}}$:
\[
\begin{array}{ll}
\underset{\mathbf{a}, \mathbf{b}, \{j_t\}\subseteq [J]}
{\mbox{\bf maximize}}  & \epc
\displaystyle{
\sum_{t \in \mathcal T_\ell}
} \, \displaystyle{
\sum_{s=1}^N
} \, \onebld(Y_s = j_t) {\underset{k \in A_R(t)}{\prod}} 
\onebld_{[ \, 0, \infty )}( 
(a^k )^\top X^s - b_k )
{\underset{k \in A_L(t)}{\prod}} \onebld_{[\, 0, \infty )}(
- (a^k )^\top X^s +  b_k 
- \varepsilon)  \\ [0.25in]
& \epc - \left( \, \lambda \, \displaystyle{
\sum_{t \in {\cal T}_{\cal B}}
} \, R_t(a^t)  \ \mbox{\begin{tabular}{l}
regularizer as a complexity measure \\
of the classification tree
\end{tabular}} \right). 
\end{array}
\] 
Note that the objective is equivalent to the minimization 
of the weighted sum of the total classification error and 
the complexity of the tree.  Moreover, given 
$( \mathbf{a},\mathbf{b} )$, the maximization over 
$\{ j_t \}_{t \in \mathcal T_\ell}$ must be achieved at the 
class that occurs most often among the data points that 
fall into the leaf node and hence it is further equivalent 
to the following program, with 
$\overline{Y}_{sj} \triangleq \left\{ \begin{array}{ll}
1 & \mbox{if $Y_s = j$} \\
0 & \mbox{otherwise}
\end{array} \right.$:
\begin{equation} \label{eq:class_tree_indicator_problem}	
\underset{\mathbf{a}, \mathbf{b}}{\mbox{\bf maximize}} 
\quad  \displaystyle{
\sum_{t \in \mathcal T_\ell}
} \, \max_{j \in [J]} \left\{ \displaystyle{
\sum_{s=1}^N
}  \,  \overline{Y}_{sj}   \onebld_{[0, \infty)} 
\left( \min\left\{\begin{array}{l} 
\underset{k \in A_R(t)}{\min} ( ( a^k )^\top X^s - b_k ), 
\\ [0.15in]     
\underset{k\in A_L(t)}{\min}(  -( a^k )^\top X^s + b_k - 
\varepsilon )
\end{array} \right\}  \right) \right\}  -  \lambda \,
\displaystyle{
\sum_{t \in {\cal T}_{\cal B}}
} \, R_t(a^t),
\end{equation}
which is thus equivalent to the MIP formulation for optimal  
tree classification in the reference \cite{BertsimasDunn17}. 
Considering the product $J \times | {\cal T}_{\ell} |$ possibly not 
too large, we could solve the above problem by rearranging the above
maximizations.  Specifically, for each $t \in {\cal T}_{\ell}$,
let $[ \, J \, ]_t$ be a copy of the index set $[ \, J \, ]$, 
and let $\mathbf{J} \triangleq \displaystyle{
\prod_{t \in {\cal T}_{\ell}}
} \, [ \, J \, ]_t$.  We write each tuple in $\mathbf{J}$ as 
$\mathbf{j} \triangleq ( j_t )_{t \in {\cal T}_{\ell}}$.  
Problem (\ref{eq:class_tree_indicator_problem}) is then equivalent
to
\[
\displaystyle{
\operatornamewithlimits{\mbox{\bf maximum}}_{
\mathbf{j}   \in \mathbf{J}}
}   \left[   \underbrace{
\underset{\mathbf{a}, \mathbf{b}}{\mbox{\bf maximize}} 
  \left\{  
\displaystyle{
\sum_{t \in \mathcal T_\ell}
} \,  \displaystyle{
\sum_{s=1}^N
}  \, \overline{Y}_{s j_t}   \onebld_{[0, \infty)} 
\left( \underbrace{\min\left\{\begin{array}{l} 
\underset{k \in A_R(t)}{\min} ( ( a^k )^\top X^s - b_k ), 
\\ [0.15in]     
\underset{k\in A_L(t)}{\min}(  -( a^k )^\top X^s + b_k - 
\varepsilon )
\end{array} \right\}}_{\mbox{concave, piecewise affine in
$(a^k,b_k)$}} 
\right)   -  \lambda  \displaystyle{
\sum_{t \in {\cal T}_{\cal B}}
} \, R_t(a^t) \, \right\}}_{\mbox{each an unconstrained 
HSCOP of type (\ref{eq:HSCOP})}} \, \right]
\]
which decomposes into $| \mathbf{J} |$ independent HSCOPs
each amenable to the theory and methods to be presented 
subsequently.  Such decomposition could facilitate the
solution of the overall problem by parallel computing.

\gap 

It is worth mentioning that we can utilize the 
decision tree as the parameterized policy class for treatment 
learning with Gini constraint that is presented in 
Section~\ref{sec:Treatment with Income}.  Under the basic setting 
of treatment learning, considering a tree-based treatment rule 
$g: \mathcal X \to [J]$ with associated branching 
parameters $\{(a^k, b^k)\}$ and assigned labels 
$\{j_t\}_{t \in \mathcal T_{\ell}} $ for all leaf nodes, we obtain
\[
\onebld(g(X) = j) = \sum_{t \in \mathcal T_{\ell}} \Bigg( \, 
\onebld(j_t = j) \, \onebld_{[0, \infty)}\left( \, 
\min \left\{ \begin{array}{l}
\displaystyle{
\min_{k \in A_R(t)}
} \, ( \, (a^k)^\top X^s - b^k \, ) , \\ [0.1in] 
\displaystyle{
\min_{k\in A_L(t)}
} \, ( \, -(a^k)^\top X^s + b^k - \varepsilon \, )  \, 
\end{array} \right\} \, \right) \, \Bigg);  
\] 
hence we can similarly obtain a HSCOP for tree-based treatment 
learning problem with the Gini constraint involving the true-false 
indicator functions $\onebld(j_t = j)$ as weights. More importantly, 
unlike  heuristic tree approaches that typically can treat only problems 
without constraints, Heaviside formulations easily  accommodate a host 
of realistic constraints, such as income  fairness in welfare economics. 

\section{Local Properties of the HSCOP: A Summary}
\label{sec:HSCOP}
\vspace{0.5ex}

Recall that for an optimization problem, stationarity conditions
are necessary for a local optimizer.  For the HSCOP (\ref{eq:HSCOP}),
a broad stationarity concept was defined in the recent
work \cite{HanCuiPang23} 
by lifting the problem to a higher-dimensional space wherein a 
standard directional stationary concept is applicable.  Specifically,
we have the following definition.

\begin{definition} \label{df:epi-stationary} \rm
A feasible vector $\bar{x} \in X_{\rm AHC}$ is an 
\emph{epi-stationary solution} of (\ref{eq:HSCOP}) if the pair
$\bar{z} \triangleq (\bar{x},f(\bar{x}))$ is a Bouligand stationary 
(abbreviated as B-stationary) solution of 
the \emph{lifted formulation}:
\begin{equation} \label{eq:lifted AHC optimization}
\displaystyle{
\operatornamewithlimits{\mbox{\bf maximize}}_{(x,t) \in {\cal Z}}
} \ \theta(x) + t 
\end{equation}
where ${\cal Z} \, \triangleq \, \left\{ \, ( x,t ) \in X_{\rm AHC} 
\times \mathbb{R} \mid t \, \leq \, f(x) \, \right\}$; i.e.,
$\theta^{\, \prime}(\bar{x};v) + dt \, \leq \, 0$ for all $( v,dt ) 
\, \in \, {\cal T}(\bar{z};{\cal Z})$,
with ${\cal T}(\bar{z};{\cal Z})$ being the tangent cone of 
${\cal Z}$ at $\bar{z}$.  Substituting the definition
of the tangent pair $( v,dt ) \, \in \, {\cal T}(\bar{z};{\cal Z})$, 
it follows that $\bar{x}$ is an epi-stationary solution 
of (\ref{eq:HSCOP}) if the implication below holds:
\begin{equation} \label{eq:epi-stat in limits}
\left. \begin{array}{l}
\displaystyle{
\lim_{\substack{k \to \infty \\ x^k \, \in \, X_{\rm AHC}}}
} \, x^k \, = \, \bar{x}; \ \displaystyle{
\lim_{k \to \infty}
} \, t_k \, = \, f(\bar{x}); \ \displaystyle{
\lim_{k \to \infty}
} \, \tau_k \, \downarrow \, 0 \\ [0.25in]
\displaystyle{
\lim_{k \to \infty}
} \, \displaystyle{
\frac{x^k - \bar{x}}{\tau_k}
} \, = \, v; \ \displaystyle{
\lim_{k \to \infty}
} \, \displaystyle{
\frac{t_k - f(\bar{x})}{\tau_k}
} \, = \, dt \\ [0.25in]
\mbox{and } \ t_k \, \leq \, f(x^k), \ \forall \, k
\end{array} \right\} \ \Rightarrow \ \theta^{\, \prime}(\bar{x};v) + 
dt \, \leq \, 0.
\end{equation}
\end{definition}

While it is possible to give an explicit description of the tangent
cone ${\cal T}(\bar{z};{\cal Z})$ 
(see \cite[Proposition~7]{CuiLiuPang22-piecewise}), a more relevant
issue for our purpose is that of \emph{stationarity sufficiency}, 
that 
is, whether there are sufficient conditions that will ensure an 
epi-stationary solution of (\ref{eq:HSCOP}) to be a local maximizer.
Rather than providing an answer to this question in its generality,
we state the following result pertaining to the problem 
(\ref{eq:HSCOP})
assuming the concavity of $\theta$, the nonnegativity of the
coefficients $\{ \psi_{ij} \}$, and the piecewise affine property
of the functions $\{ \phi_{ij} \}$.  Notice that the latter functions
are not assumed to concave or convex in this result.  The proof 
relies on the local
convexity-like property of the lifted set ${\cal Z}$ at the 
given pair $\bar{z}$; i.e., there exists a neighborhood 
${\cal N}_{\bar{z}}$ 
such that ${\cal Z} \cap {\cal N}_{\bar{z}} \subseteq \bar{z} + 
{\cal T}(\bar{z};{\cal Z})$; see \cite{HanCuiPang23}, in particular
Theorem~2 therein, for details.

\begin{proposition} \label{pr:convxity-like} \rm
Let $\theta$ be a concave function.  Suppose that the coefficients
$\left\{ \{ \psi_{ij} \}_{j=1}^{J_i} \right\}_{i=0}^m$ are nonnegative
and that the functions 
$\left\{ \{ \phi_{ij} \}_{j=1}^{J_i} \right\}_{i=0}^m$ are all 
piecewise affine.  Then $\bar{x} \in X_{\rm AHC}$ is a local 
maximizer of 
(\ref{eq:HSCOP}) if and only if $\bar{x}$ is an epi-stationary point.
\hfill $\Box$
\end{proposition}

While the concept of epi-stationarity does not have direct relevance 
to the rest of the paper, it nevertheless answers an important
question pertaining to the local optimality of a solution to
(\ref{eq:HSCOP}) from a continuous optimization perspective.  
More interestingly, the above \emph{local} result can be 
contrasted with the subsequent 
Proposition~\ref{pr:global implies fixed point} that 
pertains to a \emph{global} maximizer of (\ref{eq:HSCOP}).   We 
postpone this comparison until after the latter proposition.
Previous results related to these optimality conditions
can be found in \cite[Sections~4 and 5]{CuiLiuPang22-piecewise}.

\section{The PIP Method: Leveraging IPs of Reduced Sizes}
\label{sec:the PIP method}
\vspace{0.5ex}

Toward the goal of introducing the progressive (mixed) integer 
programming (PIP) method, we first give an 
equivalent mixed-integer programming formulation of the problem 
(\ref{eq:HSCOP}) by introducing an integer variable for each Heaviside
term.  We make the following assumption:
%

\gap

\noindent $\bullet $ there exists a constant $\underline{B}$
such that $\phi_{ij}(x) \geq \underline{B}$ for all 
$x \in P$. 
%
%

\gap

\noindent The bound $\underline{B}$ is the result
of restricting the search of the variable $x$ to 
belong to a given bounded domain, e.g.\ in the form of a bound
constraint $\| \, x \, \| \leq \overline{U}$, or
when the polyhedron $P$ itself is bounded.
Employing $\underline{B}$, we define the following mixed integer
program (MIP) involving both the continuous variable $x$ and the 
auxiliary binary variable $z$:
\begin{equation} \label{eq:all IP} 
\begin{array}{lll}
\displaystyle{
\operatornamewithlimits{\mbox{\bf maximize}}_{
x \, \in \, P; \, z}
} & \theta(x) + \displaystyle{
\sum_{j=1}^{J_0}
} \, \psi_{0j} \, z_{0j} & \\ [0.2in]
\mbox{\bf subject to} & \displaystyle{
\sum_{j=1}^{J_i}
} \, \psi_{ij} \, z_{ij} \, \geq \, b_i,
& \forall \ i \, = \, 1, \cdots, m \\ [0.25in]
& \phi_{ij}(x) \, \geq \, \underline{B} \, ( \, 1 - z_{ij} \, ), 
& \forall \ j \, = \, 1, \cdots, J_i; \ i \, = \, 0, 1, \cdots, m 
\\ [0.1in]
\mbox{\bf and} & z_{ij} \, \in \, \{ \, 0,1 \, \},
& \forall \ j \, = \, 1, \cdots, J_i; \ i \, = \, 0, 1, \cdots, m.
\end{array} \end{equation}
We note that if $\theta$ is concave and piecewise linear and so is
each function $\phi_{ij}$, then (\ref{eq:all IP}) is a mixed 
integer linear program.  When $\phi_{ij}$ is convex, the constraint
involving $\phi_{ij}$ is of the reverse convex kind.  In Appendix~A,
we present a formulation of the constraint $\phi_{ij}(x) \geq
\underline{B} ( 1 - z _{ij} )$ as a convex constraint when $\phi_{ij}$
is a piecewise affine, by introducing additional integer variables.
In this and the next section of the paper, we make no 
particular assumption on these functions unless otherwise necessary.

\gap

In contrast to the local properties in Proposition~\ref{pr:convxity-like},
the following result asserts the equivalence of the two problems
(\ref{eq:HSCOP}) and (\ref{eq:all IP}) in terms of their global 
maximizers, for arbitrary functions $\phi_{ij}$.

\begin{proposition} \label{pr:HSCOP IP form} \rm
Suppose that the coefficients
$\left\{ \{ \psi_{ij} \}_{j=1}^{J_i} \right\}_{i=0}^m$ are nonnegative. 
Let $\underline{B}$ be a lower bound of all $\phi_{ij}$ on $P$.  
If $\bar{x}$ is a global maximizer
of (\ref{eq:HSCOP}), then $( \bar{x},\bar{z} )$ is
a global maximizer of (\ref{eq:all IP}), where
$\bar{z}_{ij} \triangleq 
\onebld_{[ \, 0,\infty )}( \phi_{ij}(\bar{x}) )$ for all $(i,j)$.
Conversely, if $( \bar{x},\bar{z} )$ is
a global maximizer of (\ref{eq:all IP}), then $\bar{x}$ is a
global maximizer of (\ref{eq:HSCOP}).
\end{proposition}

\proof{Proof.}  Suppose $\bar{x}$ is an optimal solution of 
(\ref{eq:HSCOP}).   Let $(x,z)$ be feasible to (\ref{eq:all IP}).
The inequality $\phi_{ij}(x) \geq \underline{B} \, 
( \, 1 - z_{ij} \,  )$ yields that if $\phi_{ij}(x) < 0$, then 
we must have 
$z_{ij} = 0 = \onebld_{[ \, 0,\infty \, )}(\phi_{ij}(x))$.   Hence
\[
\displaystyle{
\sum_{j=1}^{J_i}
} \, \psi_{ij} \, \onebld_{[ \, 0,\infty \, )}(\phi_{ij}(x)) 
\, \geq \,  \displaystyle{
\sum_{j=1}^{J_i}
} \, \psi_{ij} \, z_{ij}, \epc \forall \, i \, = \, 0, 1, \cdots, m,
\]
with equality holding for $(\bar{x},\bar{z})$.  Thus $x$ is 
feasible to (\ref{eq:HSCOP}).  Hence
\[
\theta(\bar{x}) + \displaystyle{
\sum_{j=1}^{J_0}
} \, \psi_{0j} \, \onebld_{[ \, 0,\infty \, )}(\phi_{0j}(\bar{x})) 
\, \geq \,  \theta(x) + \displaystyle{
\sum_{j=1}^{J_0}
} \, \psi_{0j} \, \onebld_{[ \, 0,\infty \, )}(\phi_{0j}(x)).
\]
This is enough to show that  $( \bar{x},\bar{z})$ is an optimal 
solution of (\ref{eq:all IP}).

\gap

Conversely, suppose that $( \bar{x},\bar{z})$ is an optimal solution
of (\ref{eq:all IP}).  Let $x$ be feasible to (\ref{eq:HSCOP}).  
Then $(x,z)$ is feasible to (\ref{eq:all IP}), where 
$z_{ij} \triangleq \onebld_{[ \, 0,\infty \, )}(\phi_{ij}(x))$ 
for all $(i,j)$.  Hence
\[ \begin{array}{lll}
\theta(\bar{x}) + \displaystyle{
\sum_{j=1}^{J_i}
} \, \psi_{ij} \, \onebld_{[ \, 0,\infty \, )}(\phi_{ij}(\bar{x})) 
& = & \theta(\bar{x}) + \displaystyle{
\sum_{j=1}^{J_i}
} \, \psi_{ij} \, \bar{z}_{ij} \\ [0.2in]
& \geq &  \theta(x) + \displaystyle{
\sum_{j=1}^{J_i}
} \, \psi_{ij} \, z_{ij} \, = \, \theta(x) + \displaystyle{
\sum_{j=1}^{J_i}
} \, \psi_{ij} \, \onebld_{[ \, 0,\infty \, )}(\phi_{ij}(x)),
\end{array} \]
showing that $\bar{x}$ is optimal to (\ref{eq:HSCOP}).
 \endproof

\subsection{The computational workhorse}

The drawback of the formulation (\ref{eq:all IP}) is its 
computational difficulty when the number of Heaviside functions which is the same as the number of integer variables in (\ref{eq:all IP})  is 
very large.  In what follows, we introduce a restricted MIP 
formulation involving only a subset of integer variables.
Let $\varepsilon_k \geq 0$ for $k =1,2$ be given nonnegative scalars, 
presumably very small.  Given a vector $\bar{x}$ feasible to 
(\ref{eq:HSCOP}), define three complementary index sets:
\begin{equation} \label{eq:3 index sets} \begin{array}{lll}
{\cal J}_{\varepsilon_2}^< (\bar{x}) & \triangleq & \left\{ \, ( i,j )
\, \mid \, \phi_{ij}(\bar{x}) \, < \,  -\varepsilon_2 \, \right\} ,\, 
{\cal J}_{\varepsilon_{1,2}}^{\rm in}(\bar{x})   \triangleq 
\left\{ \, ( i,j ) \, \mid \,  -\varepsilon_2 \, \leq \,
\phi_{ij}(\bar{x}) \, \leq \, \varepsilon_1 \, \right\},\\[0.1in]
{\cal J}_{\varepsilon_1}^>(\bar{x}) & \triangleq & \left\{ \, ( i,j ) 
\, \mid \, \phi_{ij}(\bar{x})) \, > \, \varepsilon_1 \, \right\} .
\end{array}
\end{equation}
We term ${\cal J}_{\varepsilon_{1,2}}^{\rm in}(\bar{x})$ the 
``in-between'' index set as an enlargement of the  set $
{\cal J}_0^=(\bar{x}) \, \triangleq \, \left\{ \, ( i,j ) \, \mid \,
\phi_{ij}(\bar{x}) \, = \, 0 \, \right\}.$
The definition of these index sets is derived from the observation 
that if $\phi_{ij}(\bar{x})$ is nonzero, then for all
$x$ sufficiently near $\bar{x}$, the sign
of $\phi_{ij}(x)$ follows that of $\phi_{ij}(\bar{x})$ by continuity; 
whereas if $\phi_{ij}(\bar{x})$ is zero, then the
sign of $\phi_{ij}(x)$ cannot be determined for $x$ near $\bar{x}$. 
In practical applications, one may expect
that  the sets ${\cal J}_0^=(\bar{x})$
and also ${\cal J}_{\varepsilon_{1,2}}^{\rm in}(\bar{x})$ for 
small $\{\varepsilon_k\}_{k=1,2}$ may contain only a small 
number of elements.  Thus it is reasonable to propose
an algorithm based on solving a mixed integer subproblem at each 
iteration that involves only the integer variables $z_{ij}$ for 
$(i,j) \in {\cal J}_{\varepsilon_{1,2}}^{\rm in}(\bar{x})$ 
at an iterate $\bar{x}$.  Such a subproblem is further motivated by
the desire to improve the quality of the candidate solution $\bar{x}$
within a region in which the sign of $\phi_{ij}$ does not change 
for $(i,j) \in \mathcal J_{\varepsilon_2}^<(\bar x) \bigcup \mathcal 
J_{\varepsilon_1}^>(\bar x)$.  Therefore, we construct an IP
by assigning the integer variables only for those Heaviside 
terms with indices in 
$\mathcal J_{\varepsilon_{1,2}}^{\rm in}(\bar x)$ 
and fixing the integer variables with indices in 
${\cal J}_{\varepsilon_{1,2}}^{\rm in}(\bar{x})$ at the values 
dictated by $\bar{x}$.  This idea leads to the definition of 

\gap

\noindent {\bf The restricted integer program:}. Given the vector 
$\bar{x} \in X_{\rm AHC}$ and the three index sets 
${\cal J}_{\varepsilon_2}^< (\bar{x})$,
${\cal J}_{\varepsilon_{1,2}}^{\rm in} (\bar{x})$, and 
${\cal J}_{\varepsilon_1}^> (\bar{x})$, define the 
integer program:
\begin{equation} \label{eq:reduced IP formulation}
\begin{array}{lll}
\displaystyle{
\operatornamewithlimits{\mbox{\bf maximize}}_{x \in P; \,  z}
} & \theta(x) + \displaystyle{
\sum_{(0,j) \in {\cal J}_{\varepsilon_{1,2}}^{\, \rm in}(\bar{x})}
} \, \psi_{0j} \, z_{0j} + \displaystyle{
\sum_{(0,j) \in {\cal J}_{\varepsilon_1}^>(\bar{x})}
} \, \psi_{0j} & \\ [0.25in]
\mbox{\bf subject to} & \displaystyle{
\sum_{(i,j) \in {\cal J}_{\varepsilon_{1,2}}^{\, \rm in}(\bar{x})}
} \, \psi_{ij} \, z_{ij} + \displaystyle{
\sum_{(i,j) \in {\cal J}_{\varepsilon_1}^>(\bar{x})}
} \, \psi_{ij} \, \geq \, b_i, & i \, = \, 1, \cdots, m \\ [0.25in]
& \epc \phi_{ij}(x) \, \geq \, \underline{B} \, ( \, 1 - z_{ij} \,  ), 
& ( i,j ) \in  {\cal J}_{\varepsilon_{1,2}}^{\, \rm in}(\bar{x})
\\ [0.1in]
& \epc \phi_{ij}(x) \, \geq \, 0 \epc ( \Leftrightarrow z_{ij} = 1 ) 
& ( i,j ) \in  {\cal J}_{\varepsilon_1}^>(\bar{x}) \\ [0.1in]
& \epc \phi_{ij}(x) \ \mbox{ free} \epc ( \Leftrightarrow z_{ij} = 0 ) 
& ( i,j ) \in  {\cal J}_{\varepsilon_1}^<(\bar{x}) 
\\ [0.1in]
\mbox{\bf and} &  \epc z_{ij} \ \in \, \{ \, 0, 1 \, \}, & 
( i,j ) \in  {\cal J}_{\varepsilon_{1,2}}^{\, \rm in}(\bar{x}).
\end{array} 
\end{equation}

Similar to Proposition~\ref{pr:HSCOP IP form}, we can prove that 
(\ref{eq:reduced IP formulation}) is equivalent to the following 
partial Heavisde formulation wherein the integer variables are 
removed, leaving only the original variable $x$ in the problem:
\begin{equation} \label{eq:HSCOP reduced}
\begin{array}{lll}
\displaystyle{
\operatornamewithlimits{\mbox{\bf maximize}}_{x \in P}
} & \theta(x) + \displaystyle{
\sum_{(0,j) \, : \, \phi_{0j}(\bar{x}) \geq -\varepsilon_2}
} \, \psi_{0j} \, \onebld_{[ \, 0,\infty \, )}(\phi_{ij}(x)) & 
\\ [0.25in]
\mbox{\bf subject to} & \displaystyle{
\sum_{(i,j) \, : \, \phi_{ij}(\bar{x}) \geq -\varepsilon_2}
} \, \psi_{ij} \, \onebld_{[ \, 0,\infty \, )}(\phi_{ij}(x)) 
\, \geq \, b_i, & i \, = \, 1, \cdots, m \\ [0.25in]
\mbox{\bf and} & \epc \ \phi_{ij}(x) \, \geq \, 0, & 
( i,j ) \in  {\cal J}_{\varepsilon_1}^>(\bar{x}).
\end{array} 
\end{equation}
Let $X_{\varepsilon_{1,2}}(\bar x)$ denote the feasible region 
of (\ref{eq:HSCOP reduced}) and $\mu_{\varepsilon_{1,2}}(\bar x)$
its optimal objective value.  On one hand, it is easy to see that if 
$0 \leq (\varepsilon_1,\varepsilon_2) \leq 
(\varepsilon_1^{\, \prime},\varepsilon_2^{\, \prime})$, then
\[
\left\{ \, (i,j) \, \mid \, \phi_{ij}(\bar{x}) \, \geq \, 
-\varepsilon_2 \, \right\} \, \subseteq \,
\left\{ \, (i,j) \, \mid \, \phi_{ij}(\bar{x}) \, \geq \, 
-\varepsilon_2^{\, \prime} \, \right\}  \mbox{ and } \
{\cal J}_{\varepsilon_1^{\, \prime}}^>(\bar{x}) \, \subseteq 
{\cal J}_{\varepsilon_1}^>(\bar{x}), 
\]
which implies $X_{\varepsilon_{1,2}}(\bar x) \subseteq 
X_{\varepsilon_{1,2}^{\, \prime}}(\bar x)$, 
provided that $\psi_{ij} \geq 0$ for all $(i,j)$. 
Hence $\mu_{\varepsilon_{1,2}}(\bar x) \leq 
\mu_{\varepsilon_{1,2}^{\, \prime}}(\bar x)$.  On the other hand,
we have $\bar{x} \in X_{\varepsilon_{1,2}}(\bar x) \subseteq 
X_{\rm AHC}$ for all nonnegative pairs
$(\varepsilon_1,\varepsilon_2)$.  Thus 
$\mu_{\varepsilon_{1,2}}(\bar x) $ is a valid lower bound of the 
optimal objective value of (\ref{eq:HSCOP}).  Moreover,
since there is only a finite number of these index sets, it follows 
that there exists $\bar{\varepsilon} > 0$ such that for all
$(\varepsilon_1,\varepsilon_2) \geq 
(\bar{\varepsilon},\bar{\varepsilon})$, the corresponding optimal 
value $\mu_{\varepsilon_{1,2}}(\bar x)$ will reach the globally
optimal objective value of the Heavyside optimization problem 
(\ref{eq:HSCOP}).  The following result shows how problem
(\ref{eq:HSCOP reduced}), or equivalently
(\ref{eq:reduced IP formulation}) can provide a necessary condition
for the global optimality of a feasible solution of (\ref{eq:HSCOP});
statement (B) is particularly relevant for the PIP method to be 
presented later.

\begin{proposition} \label{pr:global implies fixed point} \rm
Assume that the functions $\phi_{ij}$ are bounded below by 
$\underline{B}$ on $P$ and that  the coefficients 
$\psi_{ij}$ are nonnegative.  Let $\bar{x}$ be a feasible solution 
of (\ref{eq:HSCOP}) and let $\bar{z}_{ij} \triangleq 
\onebld_{[ \, 0,\infty \, )}(\phi_{ij}(\bar{x}))$ for all $(i,j)$.  
Then the two statements (A) and (B) hold.

\gap

\noindent {\bf (A)} $\bar{x}$ is a globally optimal solution of 
(\ref{eq:HSCOP}) (equivalently, $( \bar{x},\bar{z} )$ 
is  a globally optimal pair of (\ref{eq:all IP})),
if and only if $\bar{x}$ is a globally optimal solution of 
(\ref{eq:HSCOP reduced})) (equivalently, 
$( \bar{x},\bar{z} )$ is a globally optimal pair of 
(\ref{eq:reduced IP formulation})) for {\bf all} 
positive $\varepsilon_1$ and $\varepsilon_2$ (sufficiently large).

\gap

\noindent {\bf (B)} $\bar{x}$ is a globally optimal solution of 
(\ref{eq:HSCOP reduced}) (equivalently, $( \bar{x},\bar{z} )$ is 
a globally optimal pair of (\ref{eq:reduced IP formulation})) for 
{\bf some} $(\varepsilon_1,\varepsilon_2) > 0$,
if and only if $\bar{x}$ is a locally optimal 
solution of (\ref{eq:HSCOP}).
\end{proposition}

\proof{Proof.} ``Only if statement in (A):''   
Suppose $\bar{x}$ is a globally optimal solution of 
(\ref{eq:HSCOP}). Let $(x,z)$ be an arbitrary feasible pair of 
(\ref{eq:reduced IP formulation}) for a given (but arbitrary) 
positive pair $(\varepsilon_1,\varepsilon_2)$.   
Since $(\bar{x},\bar{z})$ is clearly feasible to 
(\ref{eq:all IP}), and thus to (\ref{eq:reduced IP formulation}), 
it suffices to show
\begin{equation} \label{eq:objective inequality}
\displaystyle{
\sum_{(0,j) \in {\cal J}_{\varepsilon_{1,2}}^{\, \rm in}(\bar{x})}
} \, \psi_{0j} \, \bar{z}_{0j} + \displaystyle{
\sum_{(0,j) \in {\cal J}_{\varepsilon_1}^>(\bar{x})}
} \, \psi_{0j} \, \geq \, \displaystyle{
\sum_{(0,j) \in {\cal J}_{\varepsilon_{1,2}}^{\, \rm in}(\bar{x})}
} \, \psi_{0j} \, z_{0j} + \displaystyle{
\sum_{(0,j) \in {\cal J}_{\varepsilon_1}^>(\bar{x})}
} \, \psi_{0j}.
\end{equation}
For this, we verify that for all $ i = 0, 1, \cdots, m$,
\[
\displaystyle{
\sum_{j=1}^{J_i}
} \, \psi_{ij} \, \onebld_{[ \, 0,\infty \, )}( \phi_{ij}(x) ) 
\, \geq \, \displaystyle{
\sum_{(i,j) \in {\cal J}_{\varepsilon_{1,2}}^{\, \rm in}(\bar{x})}
} \, \psi_{ij} \, z_{ij} + \displaystyle{
\sum_{(i,j) \in {\cal J}_{\varepsilon_1}^>(\bar{x})}
} \, \psi_{ij};
\]
for $i = 1, \cdots, m$,
this inequality shows that $x$ is feasible to (\ref{eq:HSCOP}); 
and for $i = 0$, the inequality
establishes (\ref{eq:objective inequality}) in view of the global 
optimality of $\bar{x}$ for (\ref{eq:HSCOP}); that is 
\[ \begin{array}{lll}
\displaystyle{
\sum_{j=1}^{J_0}
} \, \psi_{0j} \, \onebld_{[ \, 0,\infty \, )}( \phi_{0j}(\bar x) ) & \geq &  \displaystyle{
\sum_{j=1}^{J_0}
} \, \psi_{0j} \, \onebld_{[ \, 0,\infty \, )}( \phi_{0j}( x) )  
\geq  \displaystyle{
\sum_{(0,j) \in {\cal J}_{\varepsilon_{1,2}}^{\, \rm in}(\bar{x})}
} \, \psi_{0j} \, z_{0j} + \displaystyle{
\sum_{(0,j) \in {\cal J}_{\varepsilon_1}^>(\bar{x})}
} \, \psi_{0j}.
\end{array}
\]
``If statement in (A):''  This is obvious because if $\varepsilon_1$ 
and $\varepsilon_2$ are sufficiently large,
then the index set ${\cal J}_{\varepsilon_{1,2}}^{\, \rm in}(\bar{x})$ 
contains all pairs $(i,j)$ for $j = 1, \cdots, J_i$
and $i = 1, \cdots m$.

\gap

\noindent ``Only if statement in (B)":  With the given pair $(\varepsilon_1, \varepsilon_2)$, let  
${\cal N}_{\varepsilon_{1,2}}(\bar x) \triangleq 
\left\{ x^{\, \prime} \in P \, \mid \, \| \bar{x}-x^{\, \prime} \, \| 
\leq \delta_{\varepsilon_{1,2}} \right\}$, where
\[ 
\delta_{\varepsilon_{1,2}} \, \triangleq \,   
\sup \left\{ \, \delta \geq 0 \, \left| \begin{array}{ll} 
\phi_{ij}(x^{\, \prime}) > 0  \ \mbox{ for all } \ 
\| \bar{x}-x^{\, \prime} \| \leq \delta, \ (i,j) \in 
{\cal J}_{\varepsilon_1}^>(\bar{x}), \mbox{ and}\\ [0.1in] 
\phi_{ij}(x^{\, \prime}) < 0 \ \mbox{ for all } \ 
\| \bar{x}-x^{\, \prime} \| \leq \delta, \ (i,j) \in 
{\cal J}_{\varepsilon_2}^<(\bar{x})
\end{array} \right. \right\}.
\]
Let $x^{\, \prime} \in {\cal N}_{\varepsilon_{1,2}}(\bar x)$ 
be feasible to (\ref{eq:HSCOP}). By the construction of $\delta_{\varepsilon_{1,2}}$, we have 
\[
\displaystyle{
\sum_{j=1}^{J_i}
} \, \psi_{ij} \, \onebld_{[ \, 0,\infty \, )}( \phi_{ij}(x') ) = \displaystyle{
\sum_{(i,j): \phi_{ij}(\bar x) \geq -\varepsilon_2}
} \, \psi_{ij} \, \onebld_{[ \, 0,\infty \, )}( \phi_{ij}(x') ),\] and thus $x^{\, \prime}$ must 
be feasible to (\ref{eq:HSCOP reduced}).  
Since $\bar{x}$ is globally optimal for the latter problem
restricted to ${\cal N}_{\varepsilon_{1,2}}(\bar x)$ which 
contains $x^{\, \prime}$, we have
\[ \begin{array}{l}
\theta(\bar{x}) + \displaystyle{
\sum_{j=1}^{J_0}
} \, \psi_{0j} \, \onebld_{[ \, 0,\infty \, )}(\phi_{ij}(\bar{x})) \\ [0.15in]
= \, \theta(\bar{x}) + \displaystyle{
\sum_{(0,j) : \phi_{ij}(\bar{x}) \geq -\varepsilon_2}
} \, \psi_{0j} \, \onebld_{[ \, 0,\infty \, )}(\phi_{ij}(\bar{x})) \\ [0.1in]
\geq \, \theta(x^{\, \prime}) + \displaystyle{
\sum_{(0,j) : \phi_{ij}(\bar{x}) \geq -\varepsilon_2}
} \, \psi_{0j} \, \onebld_{[ \, 0,\infty \, )}(\phi_{ij}(x^{\, \prime})) \, = \, 
\theta(x^{\, \prime}) + \displaystyle{
\sum_{j=1}^{J_0}
} \, \psi_{0j} \, \onebld_{[ \, 0,\infty \, )}(\phi_{ij}(x^{\, \prime})),
\end{array} \]
showing that $\bar{x}$ is a local maximizer of (\ref{eq:HSCOP}).

\gap

\noindent ``If statement in (B):" Suppose that $\bar x$ 
is a locally optimal solution of (\ref{eq:HSCOP}); i.e., there 
exists a neighborhood $\mathcal N_{\delta}(\bar x)$ with radius 
$\delta$ such that $\bar x$ is a globally optimal solution within
that neighborhood.  Then there exists $\varepsilon_1$, $\varepsilon_2$ and 
$\delta^{\prime} (\leq \delta)$ such that for all 
$x \in \mathcal N_{\delta^{\prime}}(\bar x)$, we have $\phi_{ij}(x) \geq 0, \forall (i,j) \in \mathcal J_{\varepsilon_1}^>(\bar x)$ and  $\phi_{ij}(x) < 0, \forall 
(i,j) \in \mathcal J_{\varepsilon_2}^<(\bar x)$.  By the local 
optimality of $\bar x$ corresponding to the neighborhood 
$\mathcal N_{\delta^{\prime}}(\bar x)$, we can thus deduce that 
$\bar x$ is a global optimal solution to \eqref{eq:HSCOP reduced} 
corresponding to such a pair $(\varepsilon_1,\varepsilon_2)$. 
\hfill $\Box$
\endproof

\gap

\noindent {\bf Remark.}
Posed in different spaces (the $x$-space $\mathbb{R}^n$ for 
(\ref{eq:HSCOP reduced}) versus 
the $(x,z)$-space $\mathbb{R}^n \times \{ 0,1 \}^J$ where
$J \triangleq 
\displaystyle{
\sum_{i=0}^I
} \, J_i$ for (\ref{eq:reduced IP formulation})), the two 
problems (\ref{eq:HSCOP reduced}) and
(\ref{eq:reduced IP formulation}) may not be equivalent in terms
of their local maximizers.  The reason is that the
neighborhood of the pair $(\bar{x},\bar{z})$ in the local optimality 
of the latter problem is restricted to the region
\[
\left\{ (x,\bar{z}) \mid \| x - \bar{x} \| \leq \delta, \,
\phi_{ij}(x) \geq 0 \mbox{ for all } (i,j) \mbox{ such that }
\phi_{ij}(\bar{x}) \geq 0 \right\}
\]
whose projection onto the 
$x$-space may only be a proper subset of any neighborhood of 
$\bar{x}$ defined for (\ref{eq:HSCOP reduced}), thus
invalidating the equivalence of local maximizers for
the two problems.  \hfill $\Box$

\gap

We can now compare Proposition~\ref{pr:convxity-like} and
\ref{pr:global implies fixed point}.  The former
proposition provides a necessary and sufficient
condition of a local maximizer of (\ref{eq:HSCOP}) in terms of  
epi-stationarity.  In contrast, part (a) of 
Proposition~\ref{pr:global implies fixed point} 
provides a necessary and sufficient condition (under 
assumptions therein) of a global maximier of (\ref{eq:HSCOP}) in 
terms of a fixed-point property defined by a self-induced problem
(\ref{eq:HSCOP reduced}) or its equivalent
integer program formulation (\ref{eq:reduced IP formulation}).
Part (b) of Proposition~\ref{pr:global implies fixed point} adds
to the characterization of a \emph{local} maximizer in terms of such 
an integer program. 

More can be said when $\theta = 0$ 
in the problems (\ref{eq:HSCOP}) and (\ref{eq:HSCOP reduced}).  
This special case
includes for instance the vanilla classification problem 
(\ref{eq:HSCOP standard classification}),
or the un-penalized treatment problem 
(\ref{eq:treatment optimization}) without the regularizer.

\begin{proposition} \label{pr:theta=0} \rm
Assume that $\theta = 0$, the functions $\phi_{ij}$ are bounded 
below by $\underline{B}$ on $P$, and the coefficients 
$\psi_{ij}$ are nonnegative.  The following three statements
are equivalent for any pair $( \bar{x},\bar{z} )$, where
$\bar{z}_{ij} \triangleq 
\onebld_{[ \, 0,\infty \, )}(\phi_{ij}(\bar{x}))$ for all $(i,j)$:

\gap

\noindent {\bf (A)} $\bar{x}$ is a globally optimal solution of 
(\ref{eq:HSCOP reduced}) (equivalently, $( \bar{x},\bar{z} )$ is 
a globally optimal pair of (\ref{eq:reduced IP formulation})) for 
some $(\varepsilon_1,\varepsilon_2) > 0$;

\gap

\noindent {\bf (B)} $\bar{x}$ is a locally optimal solution of 
(\ref{eq:HSCOP});

\gap

\noindent {\bf (C)} $\bar{x}$ is a feasible solution of (\ref{eq:HSCOP}).
\end{proposition}

\proof{Proof.}  It suffices to prove (C) $\Rightarrow$ (A).
With $\bar{x}$ feasible to (\ref{eq:HSCOP}), let
$\underline{\varepsilon}_2$ be defined by
\[
0 > -\underline{\varepsilon}_2 \, \triangleq  
\mbox{\bf maximum}\left\{ \, \phi_{ij}(\bar{x}) \mid \phi_{ij}(\bar{x}) 
< 0 \, \right\}.
\]
It then follows that for every $i\in [m]$, for all 
$0 < \varepsilon_2  < \underline{\varepsilon}_2 $, 
\[
\{ j \mid \phi_{ij}(\bar{x}) \geq -\varepsilon_2 \} \, = \, \{  j \mid 
\phi_{ij}(\bar{x}) \geq 0 \}. 
\]
Thus \eqref{eq:HSCOP reduced} is equivalent to the following program,
\[
\begin{array}{lll}
\displaystyle{
\operatornamewithlimits{\mbox{\bf maximize}}_{x \in P}
} & \displaystyle{
\sum_{j \, : \phi_{0j}(\bar{x}) \geq 0} 
} \, \psi_{0j} \, \onebld_{[ \, 0,\infty \, )}(\phi_{ij}(x))    
& \\ [0.25in]
\mbox{\bf subject to} & \displaystyle{
\sum_{j \, : \,  \phi_{ij}(\bar{x}) \geq 0}
} \, \psi_{ij} \, \onebld_{[ \, 0,\infty \, )}(\phi_{ij}(x))   
\, \geq \, b_i, & i \, = \, 1, \cdots, m \\ [0.25in]
\mbox{\bf and} & \epc \ \phi_{ij}(x) \, \geq \, 0, \epc \forall \,  
( i,j ) \in  {\cal J}_{\varepsilon_1}^>(\bar{x}).
\end{array} 
\]
Clearly, $\bar x$ is feasible to the above program and also achieves 
the maximum objective value
by the nonnegativity of coefficients $\{\psi_{0j}\}$. 
\endproof

\subsection{Description the PIP method} \label{subsec:description PIP}

We present an algorithm for solving the HSCOP (\ref{eq:HSCOP}).  The 
algorithm starts with a feasible solution of the problem and
iteratively solves a sequence of mixed-integer programs 
(\ref{eq:reduced IP formulation}) with adaptive control of the 
parameter 
pairs $\{ \varepsilon^{\nu}_1,\varepsilon^{\nu}_2 )$.   We call this
a progressive IP method because the number of integer variables is
progressively increased until the objective values do not improve. 
We should point out that the description below pertains to a 
conceptual framework of the method; variations of the steps 
are possible. In particular, the change of the index sets
and the termination rule used in the numerical experiments 
in Section~\ref{sec:numerics} are slightly different from 
those described below; see the details in the later section.


\noindent\makebox[\linewidth]{\rule{\textwidth}{0.5pt}}

\noindent {\bf Initialization:} Let $x^0 \in X_{\rm AHC}$ and 
let $\mu_0$ be the corresponding objective value of (\ref{eq:HSCOP}).
Let the parameter pair $(\varepsilon_1^0, \varepsilon_2^0)$ be
arbitrary.

\gap 

\noindent  {\bf For $\nu=0, 1, \ldots,$}  

\begin{enumerate}
\item Determine ${\cal J}_{\nu}^> \triangleq 
{\cal J}_{\varepsilon_1^{\nu}}^>(x^{\nu})$ and
${\cal J}_{\nu}^{\rm in} \triangleq 
{\cal J}_{\varepsilon_{1,2}^{\nu}}^{\rm in}(x^{\nu})$. 

\item Solve the restricted MIP (\ref{eq:reduced IP formulation}) with the 
reference point $x^\nu$ to obtain a globally optimal solution $x^{\nu+1}$ 
of (\ref{eq:HSCOP reduced}) with the corresponding optimal objective 
value $\mu_{\nu+1}$.   

\item (Enlargement of the in-between index set)  If 
$\mu_{\nu+1} = \mu_{\nu}$,  let
$\varepsilon^{\nu+1}_1 > \varepsilon^{\nu}_1$
and $\varepsilon^{\nu+1}_2 > \varepsilon^{\nu}_2$; in this case, 
${\cal J}_{\nu}^{\rm in} \subseteq {\cal J}_{\nu+1}^{\rm in}$ and 
${\cal J}_{\nu+1}^> \subseteq {\cal J}_{\nu}^>$.

\item (Shrinkage of the in-between index set) If $\mu_{\nu+1} > 
\mu_{\nu}$, let $\varepsilon^{\nu+1}_1 < \varepsilon^{\nu}_1$
and $\varepsilon^{\nu+1}_2 < \varepsilon^{\nu}_2$; in this case, 
${\cal J}_{\nu+1}^{\rm in} \subseteq {\cal J}_{\nu}^{\rm in}$.
and ${\cal J}_{\nu+1}^ > \supseteq {\cal J}_{\nu}^>$.

\item Terminate the iterations if the optimal objective values
$\mu_{\nu+1}$ remain unchanged after a few consecutive expansions
of the in-between intervals. 
\end{enumerate}

\noindent {\bf end for}.

\noindent\makebox[\linewidth]{\rule{\textwidth}{0.5pt}}

\gap

\noindent Starting with a feasible solution $x^0$, 
the iterations ensure that each iterate $x^{\nu+1}$ must remain 
feasible to the original problem (\ref{eq:HSCOP}).  In general,
the initial feasible solution $x^0$ may not be readily available;
but the method can be applied to a penalization formulation 
similar to 
(\ref{eq:treatment optimization}).  Even with such a feasible
formulation, an IP solver may not be able to find a feasible solution
to the original problem 
after a long computational time.  This was witnessed by a treatment
learning problem with 1000 data points, the state-of-the-art Gurobi
IP solver was unable to find a feasible solution (with $\gamma = 0$)
within 2 hours. 
In general, with an initial feasible solution to (\ref{eq:HSCOP})
available, the PIP
method is well defined and generates a sequence of feasible solutions
with non-decreasing
objective values all of which are valid lower bounds of the maximum
objective value of the original HSCOP.  This is accomplished 
primarily due to the expansion step.  In contrast, the purpose of 
the shrinkage
step is to reduce the number of integer variables to help reduce
the computational effort of the iterations.  We have the following
important property of the iterate at termination of the algorithm,
which is an immediate consequence of part (B) of 
Proposition~\ref{pr:global implies fixed point}.

\begin{proposition} \label{pr:termination yields local max} \rm
If $\mu_{\nu+1} = \mu_{\nu}$, then $x^{\, \nu}$ is a local maximizer 
of the original Heavyside optimization problem (\ref{eq:HSCOP}). 
\hfill $\Box$
\end{proposition}

The distinct advantage of the PIP method is twofold: one, given a
feasible $x^{\nu}$ of (\ref{eq:HSCOP}), it has the potential of
improving its objective value; if there is no improvement, then
$x^{\nu}$ is a local maximizer of the problem.  Practically,
as demonstrated in the experimental results on the regularized
treatment problem (\ref{eq:treatment optimization}) with 
$\ell_1$-penalizations to be reported in the next
section, the PIP method performs favorably well compared to the
solution of the full IP formulation (\ref{eq:all IP}) by the 
Gurobi solver in two respects: except for
problems with a relatively small number (100) of Heaviside terms 
which Gurobi can solve very efficiently, for problems
with modest number (300) of integer variables,
PIP solves the problems to within small fractions of global optimality
in substantially less times than required by Gurobi solving the
full IP in one shot; for problems
with 500 integer variable, PIP already outperforms Gurobi in both 
objective values and computational times; for problems
with large number (1,000) of integer variables and smaller Gini 
coefficients, Gurobi simply cannot
find a feasible integer solution satisfying the Gini constraint 
(as demonstrated by a positive penality value $\gamma$) within 
a preset of one hour limit while PIP obtains a local maximizer
(by Proposition~\ref{pr:termination yields local max})
in 40 minutes at the most.

\section{Numerical Experiments} \label{sec:numerics}

In this section, we illustrate the empirical performance of the PIP 
algorithm on the treatment optimization 
problem (\ref{eq:treatment optimization}) with synthetic data.  The
results provides evidence and details to support the advantages of
PIP as summarized at the end of the last section.


\subsection{Setup}

Our synthetic data comprises 30 continuous covariates 
$X \in\mathcal{X}=[0,1]^{30}$ and 4 treatment arms 
$D\in\{1,2,3,4\}$. The outcomes are generated according to 
\[ \begin{array}{ll}
Y \, = & X_{5} + \exp \left\{2+0.2X_0-0.1X_1+2X_0X_1 \right\} + 
(-0.8+1.8X_1-0.2X_2) \, \mathbf{1}\{D=1\} \\ [0.1in]
& + \, (-1 + 2.1X_1-1.2X_0) \, \mathbf{1}\{D=2\} + 
(-0.8+1.3X_0X_2) \, \mathbf{1}\{D=3\} \\ [0.1in]
& + \,(-0.4+1.8X_0-1.2X_1X_2) \, \mathbf{1}\{D=4\} + \varepsilon,
\end{array} \]
with the random error $\varepsilon$ being drawn from 
$\text{Lognormal}(0,0.001)$.
We systematically vary the number of distinct combinations of 
covariates to be $25$, $75$, $125$, $250$, and $375$, correspondingly 
generating datasets of sizes $1000$, $3000$, $5000$, $10000$, and 
$15000$. Considering the number of treatment arms as $4$, the 
resulting number of indicator functions in these scenarios is 
$100$, $300$, $500$, 
$1000$, and $1500$, respectively. For each scenario, we generate 
$5$ datasets.
For each dataset, we examine two different Gini thresholds: 
$\alpha=0.7$ and $\alpha = 0.5$, with the latter inducing a higher
level of equality and, therefore, potentially resulting in a smaller
set of feasible solutions. We set $\lambda = 0.01$ and $\rho = 10^8$. We select $\underline{B}$ based on 
$\underline{B}\leq \inf h_j(X^s,\boldsymbol{\beta})$, 
where $X^s\in[0,1]^{30}$, and assuming that the parameter space 
for each $\beta^j$ to be $[-1, 1]^{30}$. Regarding the 
classification functions, we employ the shifted score functions 
$h_j^{\tau}$ in (\ref{eq:h hat}) with each
$\tau_j=0.001$.
We use the MIP (\ref{eq:all IP}) with the full set of integer 
variables as the benchmark, which we label as ``full MIP''.  
In the case of PIP, we restrict the maximum number of 
integers considered in each subproblem to be $40\%$ and $60\%$ 
of the integers in the entire problem.  If the number of integers 
exceeds the threshold, we shrink the $\varepsilon$ in-between 
interval until the specified percentage is met. For the stopping
rule of PIP, we experimented with some preliminary runs limiting
the number of in-between interval expansions to be 15.  We observed
that in these runs, the objective values tended to remain unchanged
after 10 such expansions.  We then adopted the latter number as
the stopping rule for the results that we report in 
Subsection~\ref{subsec:results}; the savings in computational time 
with 10 expansions versus 15 is significant.

We conduct our experiments in Python and solve the problems using 
Gurobi 10.0 \cite{gurobi2023}. All the problems are solved on a
device with 13th Gen Intel Core processor and 16GB of memory. We set 
a time limit of 10 minutes for problems with 100 integers, 30 minutes 
for problems with 300 integers, and 60 minutes for problems with the 
number of integers equal or larger than 500.

\subsection{Results} \label{subsec:results}

We summarize our experimental results in five tables. Our 
comparison between PIP and full MIP involves in-sample evaluations, 
where we compute the welfare and Gini coefficient achieved by the 
computed policy under the two choices of the PIP percentages. 
The abbreviation "infeas." in the tables indicates that the obtained
solution fails to meet the Gini constraint ($\gamma > 0$ 
in (\ref{eq:treatment optimization})).

The five datasets in each scenario represent a diverse set of cases. 
As depicted in Table~\ref{tab::100}, when the number of integers is 
set to 100, full MIP performs pretty good, identifying the optimal 
solution in all datasets.  In comparison to full MIP, the 
performance of PIP is less satisfactory, especially 
the 40\% setting, as it struggles to find feasible solutions in 
most instances within 10 
in-between interval expansions.  Nevertheless, the 60\% PIP 
demonstrates improved performance within very short time too, albeit 
failing to achieve (but being capable to obtain within 5\% of) 
the optimal value.

\begin{table}[!ht]
\small 
    \centering
    \caption{Integer = 100}
    \label{tab::100}
    \begin{tabular}{|c|c|ccc|ccc|}
    \hline
         &  & & \textcolor{blue}{Gini = 0.7}  &  & & \textcolor{blue}{Gini = 0.5}  &  \\[0.02in]
         \hline
        ~ & ~ & Welfare in (\ref{eq:treatment optimization}) & Gini & Time (secs) & Welfare in (\ref{eq:treatment optimization}) & Gini & Time (secs)\\[0.02in] \hline
        ~ & full MIP & 14.263 & 0.398 & 0.889 & 14.263 & 0.416 
        & 0.824 \\ [0.02in]
        Dataset 1 & PIP (0.4) & infeas.\ & infeas.\ & 0.508 & infeas.\ & infeas.\ & 0.513 \\[0.02in] 
        ~ & PIP (0.6) & 13.389 & 0.700 & 1.815 & 13.389 & 0.500 & 1.763 \\[0.02in] \hline
        ~ & full MIP & 15.933 & 0.462 & 0.783 & 15.933 & 0.500 
        & 1.763 \\ [0.02in]
        Dataset 2 & PIP (0.4) & infeas.\ & infeas.\ & 0.514 
        & infeas.\ & infeas.\ & 0.521 \\[0.02in] 
        ~ & PIP (0.6) & 15.933 & 0.700 & 2.585 & 15.933 & 0.465 
        & 2.682 \\[0.02in] \hline
        ~ & full MIP & 18.846 & 0.508 & 0.763 & 18.196 & 0.496 
        & 391.195 \\ [0.02in]
        Dataset 3 & PIP (0.4) & 18.456 & 0.700 & 1.396 & infeas.\ 
        & infeas.\ & 30.902 \\ [0.02in]
        ~ & PIP (0.6) & 18.456 & 0.700 & 1.398 & infeas.\ & infeas.\ 
        & 214.294 \\[0.02in] \hline
        ~ & full MIP & 11.497 & 0.261 & 0.784 & 11.497 & 0.261 
        & 0.763 \\ [0.02in]
        Dataset 4 & PIP (0.4) & infeas.\ & infeas.\ & 0.538 & infeas.\ 
        & infeas.\ & 0.540 \\ [0.02in]
        ~ & PIP (0.6) & 10.890 & 0.700 & 1.358 & 10.796 & 0.293 
        & 1.579 \\[0.02in] \hline
        ~ & full MIP &16.751 & 0.700 & 0.742 & 16.751 & 0.397 & 0.731 \\[0.02in] 
        Dataset 5 & PIP (0.4) & 15.595 & 0.603 & 1.141 & 15.619 & 0.500 & 1.102 \\ [0.02in]
        ~ & PIP (0.6) & 15.595 & 0.603 & 1.139 & 15.619 & 0.500 & 1.106 \\[0.02in] \hline
    \end{tabular}
\end{table}
Table~\ref{tab::300} illustrates the performance of full MIP and PIP 
with 300 integers.  Here, we are able to observe the advantage 
of PIP in terms of time, coupled with welfare outcomes that are
comparable to those obtained from full MIP, using noticeably less 
time under the 0.7 Gini constraint.  Concerning the 0.5 Gini 
constraint, while 60\% PIP achieves more satisfactory welfare, 
40\% PIP demonstrates a more apparent advantage in terms of 
computational time.  In 4 out of 5 datasets, 40\% PIP takes 
only 3\%-4\% of the time used by full MIP and achieves at least 99\% 
of the welfare obtained by full MIP. The time difference between 
Gini=0.7 and Gini=0.5 illustrates an escalation in the level of 
problem difficulty as the constraint becomes more stringent.

As the number of integers increases to 500, as shown in 
Table~\ref{tab::500}, the benefits of PIP become more pronounced, 
particularly with 40\% PIP.  Under the 0.7 Gini constraint, the time
taken by 40\% PIP is 3\%-10\% of that of full MIP, and the welfare 
achieved by PIP exceeds or is at least not smaller than the
welfare of full MIP by 1\%. 
Under the 0.5 Gini constraint, PIP achieves higher welfare in 
much less time compared to full MIP, even when the latter is feasible. 
Full MIP 
is unable to obtain a solution that satisfies the Gini constraint 
within 60 minutes in 4 out of 5 datasets.
\begin{table}[!ht]
\small 
    \centering
    \caption{Integer = 300}
    \label{tab::300}
    \begin{tabular}{|c|c|ccc|ccc|}
    \hline
     &  & & \textcolor{blue}{Gini = 0.7}  &  & & 
     \textcolor{blue}{Gini = 0.5}  &  \\[0.02in] 
     \hline
        ~ & ~ & Welfare in (\ref{eq:treatment optimization}) & Gini & Time (secs) & Welfare in (\ref{eq:treatment optimization}) & Gini & Time (secs)\\[0.02in] \hline
        ~ & full MIP &17.208 & 0.556 & 620.750 & 14.835 & 0.500 
        & 1802.784 \\[0.02in] 
        Dataset 1 & PIP (0.4) & 17.144 & 0.565 & 22.061 & 13.745 
        & 0.499 & 1216.713 \\[0.02in] 
        ~ & PIP (0.6) & 17.167 & 0.561 & 71.927 & 14.801 & 0.498 
        & 1949.680 \\[0.02in] 
        \hline
        ~ & full MIP &14.900 & 0.699 & 880.472 & 11.891 & 0.499 
        & 1802.743 \\ [0.02in]
        Dataset 2 & PIP (0.4) & 14.872 & 0.550 & 40.645 & 12.340 
        & 0.499 & 1112.538 \\[0.02in] 
        ~ & PIP (0.6) & 14.881 & 0.549 & 88.034 & 13.069 & 0.495 
        & 1790.086 \\[0.02in] 
        \hline
        ~ & full MIP &14.691 & 0.700 & 130.776 & 12.478 & 0.494 
        & 1802.741 \\[0.02in] 
        Dataset 3 & PIP (0.4) & 14.290 & 0.578 & 35.374 & 12.327 
        & 0.499 & 1010.395 \\[0.02in] 
        ~ & PIP (0.6) & 14.571 & 0.565 & 75.543 & 12.583 & 0.499 
        & 1450.673 \\ [0.02in]
        \hline
        ~ & full MIP &17.436 & 0.621 & 1802.923 & 14.341 & 0.491 
        & 1802.953 \\[0.02in] 
        Dataset 4 & PIP (0.4) & 17.350 & 0.567 & 29.542 & 14.452 
        & 0.499 & 1408.018 \\[0.02in] 
        ~ & PIP (0.6) & 17.394 & 0.564 & 136.028 & 14.496 & 0.484 
        & 1942.064 \\ [0.02in]
        \hline
        ~ & full MIP &13.085 & 0.474 & 40.612 & 13.085 & 0.498 
        & 68.329 \\[0.02in] 
        Dataset 5 & PIP (0.4) & 13.080 & 0.454 & 19.499 & 13.059 
        & 0.456 & 17.914 \\[0.02in] 
        ~ & PIP (0.6) & 13.085 & 0.453 & 31.185 & 13.085 & 0.453 
        & 37.874 \\[0.02in] \hline

    \end{tabular}
\end{table}

\begin{table}[!ht]
\small 
    \centering
    \caption{Integer = 500}
    \label{tab::500}
    \begin{tabular}{|c|c|ccc|ccc|}
        \hline
         &  & & \textcolor{blue}{Gini = 0.7}   &  & & 
         \textcolor{blue}{Gini = 0.5}  &  \\[0.02in] 
         \hline
        ~ & ~ & Welfare in (\ref{eq:treatment optimization}) & Gini & Time (secs) & Welfare in (\ref{eq:treatment optimization}) & Gini & Time (secs)\\[0.02in]  \hline
        ~ & full MIP &20.597 & 0.700 & 3608.139 & infeas.\ & infeas.\ & 3607.859 \\[0.02in]  
        Dataset 1 & PIP (0.4) & 20.614 & 0.580 & 182.499 & 16.084 & 0.500 & 948.041 \\[0.02in]  
        ~ & PIP (0.6) & 20.614 & 0.580 & 417.164 & 16.084 & 0.500 & 948.436 \\[0.02in]  \hline
        ~ & full MIP &17.019 & 0.700 & 3607.516 & infeas.\ & infeas.\ & 3607.856 \\[0.02in]  
        Dataset 2 & PIP (0.4) & 17.024 & 0.515 & 103.761 & 16.465 & 0.500 & 1155.054 \\[0.02in]  
        ~ & PIP (0.6) & 17.043 & 0.510 & 855.483 & 16.319 & 0.499 & 1364.311 \\[0.02in]  \hline
        ~ & full MIP &16.555 & 0.700 & 3607.713 & infeas.\ & infeas.\ & 3608.004 \\[0.02in]  
        Dataset 3 & PIP (0.4) & 16.475 & 0.696 & 69.403 & 15.657 & 0.497 & 1131.078 \\[0.02in]  
        ~ & PIP (0.6) & 16.475 & 0.696 & 69.325 & 15.657 & 0.497 & 1133.224 \\[0.02in]  \hline
        ~ & full MIP &16.897 & 0.697 & 3607.927 & 16.985 & 0.500 & 3608.132 \\[0.02in]  
        Dataset 4 & PIP (0.4) & 16.997 & 0.500 & 262.780 & 16.997 & 0.500 & 927.229 \\[0.02in]  
        ~ & PIP (0.6) & 16.981 & 0.505 & 1331.195 & 16.997 & 0.500 & 1045.514 \\[0.02in]  \hline
        ~ & full MIP &18.103 & 0.687 & 3608.471 & infeas.\ & infeas.\ & 3607.892 \\[0.02in]  
        Dataset 5 & PIP (0.4) & 18.103 & 0.603 & 343.598 & 13.418 & 0.500 & 1943.345 \\[0.02in]  
        ~ & PIP (0.6) & 18.103 & 0.603 & 1011.516 & 13.283 & 0.499 & 2002.627 \\[0.02in]  \hline
    \end{tabular}
\end{table}

With a further increase in the number of integers to 1000 
(Table~\ref{tab::1000}), PIP is able to achieve 
higher welfare than full MIP in significantly shorter times, ranging 
from 17 to 40 minutes, even under a more relaxed 0.7 constraint.
Under the 
stricter 0.5 constraint, full MIP fails to find a feasible solution 
within the 60-minute limit for 4 out of 5 datasets; this is 
in contrast to PIP's success for these problems.  
In Table~\ref{tab::1500}, Dataset~2 
illustrates a scenario where there might not be a feasible solution 
under the 0.5 Gini constraint.  Except for this problem, this table 
further accentuates the 
limitations of full MIP for this type of constrained optimization 
problem in finding feasible solutions when the number of integers 
is relatively large, underscoring the substantial advantage of PIP 
and its superiority with the increase in the number of integer 
variables.

To conduct a comprehensive comparison between PIP and full MIP, 
we further present the objective values of incumbent feasible
solution attained by full MIP at various time points. These time points
are chosen to be comparable to those of PIP (0.4) and PIP (0.6),
offering a basis for a detailed comparison with PIP.  In cases where 
full MIP exhibits little improvement or lacks a suitable time point 
for comparison with PIP, we document the time at which full MIP
identifies a first feasible solution.  As illustrated in 
Table~\ref{tab::time}, for instances where the number of integers 
is set at 300 and 500, the objective value achieved by full MIP
within the time frame utilized by PIP is in general lower than that 
of PIP.  With an increase in the number of integers to 1000, it 
becomes apparent that full MIP struggles to exhibit significant 
improvement over time following the identification of a feasible 
solution. If there is no improvement within the time limit of 3600 
seconds for full MIP, we display the objective value when the first 
feasible solution is obtained, as well as when the time limit 
is reached.

\gap

In summary, the advantage of PIP starts with the problems with
300 integer variables and becomes quite pronounced on the large
problems.  Specifically, except for the first set of results 
with 100 integer variables, PIP demonstrates a clear 
advantage on problems of all other sizes with the 0.7 Gini 
threshold in terms of computational times, producing 
solutions that are quite comparable (in the case of 300 integer 
variables) and superior (in the remaining cases) to the full MIP 
solutions. Even under the stricter 0.5 Gini constraint, while 
the times taken by PIP increase, it continues to produce feasible 
solutions in 
cases where full MIP fails, particularly when dealing with a large 
number of integers.  Overall, these numerical results provide
compelling evidence demonstrating the appeal of the very simple 
yet successful progressive integer approach to deal with 
otherwise unsolvable large integer programs.  Our further research
will apply the method to many nonconvex nondifferentiable
optimization problems whose state-of-the-art solution methods
essentially stop with a stationary solution of some kind.  By 
leveraging the efficiency of integer programming methods for 
problems of modest sizes, we are able to not only obtain 
high-quality improvements of such stationary solutions and provide
certificates of the global optimality of these improved solutions
in the verifiable cases.

\begin{table}[!ht]
\small 
    \centering
    \caption{Integer = 1000}
    \label{tab::1000}
    \begin{tabular}{|c|c|ccc|ccc|}
    \hline
         &  & & \textcolor{blue}{Gini = 0.7}  &  & & 
         \textcolor{blue}{Gini = 0.5}  &  \\[0.02in]
        \hline
        ~ & ~ & Welfare in (\ref{eq:treatment optimization}) & Gini & Time (secs) & Welfare in (\ref{eq:treatment optimization}) & Gini & Time (secs)\\[0.02in] \hline
        ~ & full MIP &15.302 & 0.597 & 3630.367 & infeas.\ & infeas.\ & 3631.948 \\[0.02in] 
        Dataset 1 & PIP (0.4) & 16.162 & 0.570 & 1768.071 & 13.438 & 0.500 & 2088.512 \\[0.02in] 
        ~ & PIP (0.6) & 16.162 & 0.570 & 1795.900 & 13.438 & 0.500 & 2111.811 \\[0.02in] \hline
        ~ & full MIP &14.148 & 0.699 & 3632.826 & infeas.\ & infeas.\ & 3632.248 \\ [0.02in] 
        Dataset 2 & PIP (0.4) & 14.932 & 0.593 & 1189.601 & 11.482 & 0.500 & 2310.535 \\ [0.02in] 
        ~ & PIP (0.6) & 14.932 & 0.593 & 1197.255 & 11.482 & 0.500 & 2329.943 \\[0.02in]  \hline
        ~ & full MIP &14.860 & 0.585 & 3630.763 & infeas.\ & infeas.\ & 3630.753 \\[0.02in]  
        Dataset 3 & PIP (0.4) & 15.765 & 0.554 & 1373.404 & 13.825 & 0.499 & 1371.878 \\[0.02in]  
        ~ & PIP (0.6) & 15.765 & 0.553 & 1386.319 & 13.825 & 0.499 & 1380.434 \\[0.02in]  \hline
        ~ & full MIP &14.029 & 0.570 & 3630.715 & infeas.\ & infeas.\ & 3630.745 \\ [0.02in] 
        Dataset 4 & PIP (0.4) & 15.019 & 0.531 & 1573.635 & 13.525 & 0.499 & 1981.332 \\ [0.02in] 
        ~ & PIP (0.6) & 14.992 & 0.533 & 1409.802 & 13.525 & 0.499 & 2009.177 \\[0.02in]  \hline
        ~ & full MIP &13.734 & 0.546 & 3631.306 & infeas.\ & 0.495 & 3633.033 \\ [0.02in] 
        Dataset 5 & PIP (0.4) & 14.859 & 0.511 & 1619.978 & 14.653 & 0.500 & 1989.155 \\ [0.02in] 
        ~ & PIP (0.6) & 14.859 & 0.511 & 1649.590 & 14.653 & 0.500 & 2009.995 \\[0.02in]  \hline
    \end{tabular}
\end{table}

\begin{table}[!ht]
\small 
    \centering
    \caption{Integer = 1500}
    \label{tab::1500}
    \begin{tabular}{|c|c|ccc|ccc|}
    \hline
         &  & & \textcolor{blue}{Gini = 0.7}  &  & & 
         \textcolor{blue}{Gini = 0.5}  &  \\ [0.02in] 
         \hline
        ~ & ~ & Welfare in (\ref{eq:treatment optimization}) & Gini & Time (secs) & Welfare in (\ref{eq:treatment optimization}) & Gini & Time (secs)\\[0.02in]  \hline
        ~ & full MIP &infeas.\ & infeas.\ & 3668.845 & infeas.\ & infeas.\ & 3670.011 \\[0.02in]  
        Dataset 1 & PIP (0.4) & 15.589 & 0.587 & 1203.570 & 13.041 & 0.499 & 3453.371 \\ [0.02in] 
        ~ & PIP (0.6) & 15.589 & 0.587 & 1035.136 & 13.041 & 0.499 & 3278.364 \\[0.02in]  \hline
        ~ & full MIP &15.187 & 0.631 & 3670.815 & infeas.\ & infeas.\ & 3669.271 \\ [0.02in] 
        Dataset 2 & PIP (0.4) & 16.099 & 0.616 & 2127.023 & infeas.\ & infeas.\ & 2542.495 \\ [0.02in] 
        ~ & PIP (0.6) & 16.099 & 0.616 & 2171.334 & infeas.\ & infeas.\ & 2610.988 \\[0.02in]  \hline 
        ~ & full MIP &17.897 & 0.650 & 3669.292 & infeas.\ & infeas.\ & 3669.332 \\[0.02in]  
        Dataset 3 & PIP (0.4) & 18.864 & 0.631 & 2035.662 & 13.029 & 0.500 & 2933.194 \\[0.02in]  
        ~ & PIP (0.6) & 18.864 & 0.631 & 2101.903 & 13.029 & 0.500 & 2990.741 \\[0.02in]  \hline
        ~ & full MIP &16.953 & 0.634 & 3668.571 & infeas.\ & infeas.\ & 3667.437 \\[0.02in]  
        Dataset 4 & PIP (0.4) & 18.002 & 0.612 & 2179.204 & 13.454 & 0.500 & 2488.458 \\[0.02in]  
        ~ & PIP (0.6) & 18.002 & 0.612 & 2233.348 & 13.454 & 0.500 & 2549.603 \\[0.02in]  \hline
        ~ & full MIP &13.507 & 0.583 & 3667.660 & infeas.\ & infeas.\ & 3666.726 \\[0.02in]  
        Dataset 5 & PIP (0.4) & 14.472 & 0.552 & 1852.692 & 12.471 & 0.499 & 869.962 \\[0.02in]  
        ~ & PIP (0.6) & 14.472 & 0.552 & 1918.453 & 12.471 & 0.499 & 870.097 \\[0.02in]  \hline
    \end{tabular}
\end{table}

\begin{table}[!ht]
\small 
    \centering
    \caption{Comparison between MIP and PIP based on the
    PIP times (Gini=0.7)}
    \label{tab::time}
    \begin{tabular}{|c|c|cc|cc|cc|}
    \hline
         &  & \multicolumn{2}{|c|}{\textcolor{blue}{Integer = 300}}  & \multicolumn{2}{|c|}{\textcolor{blue}{Integer = 500}} & \multicolumn{2}{|c|}{\textcolor{blue}{Integer = 1000}}   \\[0.02in]  \hline
         &  & Welfare  & Time  & Welfare  & Time  
         & Welfare  & Time   \\ [0.02in]  
         &  & in (\ref{eq:treatment optimization}) & (secs)  
         & in (\ref{eq:treatment optimization}) & (secs)  
         & in (\ref{eq:treatment optimization}) & (secs)   \\[0.02in]  \hline
         & full MIP &15.618 & 26.000 & 20.472 & 186.000 & 15.302 & 49.000  \\[0.02in]  
        Dataset 1 & full MIP &16.937 & 75.000 & 20.541 & 430.000 & 15.302 & 3600.000  \\[0.02in]  
        ~ & PIP (0.4) & 17.144 & 22.061 & 20.614 & 182.499 & 16.162 & 1768.071  \\[0.02in]  
        ~ & PIP (0.6) & 17.167 & 71.927 & 20.614 & 417.164 & 16.162 & 1795.900  \\[0.02in]  \hline
         & full MIP &13.184 & 3.000 & 14.728 & 8.000 & 13.932 & 41.000  \\[0.02in]  
        Dataset 2 & full MIP &13.184 & 96.000 & 15.508 & 1131.000 & 14.001 & 1215.000  \\[0.02in]  
        ~ & PIP (0.4) & 14.872 & 40.645 & 17.024 & 103.761 & 14.932 & 1189.601  \\[0.02in]  
        ~ & PIP (0.6) & 14.881 & 88.034 & 17.043 & 855.483 & 14.932 & 1197.255  \\[0.02in]  \hline
         & full MIP &13.958 & 35.000 & 14.562 & 8.000 & 14.860 & 44.000  \\[0.02in]  
        Dataset 3 & full MIP &14.674 & 77.000 & 14.562 & 69.000 & 14.860 & 3600.000  \\[0.02in]  
        ~ & PIP (0.4) & 14.290 & 35.374 & 16.475 & 69.403 & 15.765 & 1373.404  \\[0.02in]  
        ~ & PIP (0.6) & 14.571 & 75.543 & 16.475 & 69.325 & 15.765 & 1386.319  \\[0.02in]  \hline
         & full MIP &15.875 & 21.000 & 15.366 & 202.000 & 14.029 & 48.000  \\[0.02in] 
        Dataset 4 & full MIP &15.918 & 180.000 & 15.694 & 1359.000 & 14.029 & 3600.000  \\[0.02in]  
        ~ & PIP (0.4) & 17.350 & 29.542 & 16.997 & 262.780 & 15.019 & 1573.635  \\[0.02in]  
        ~ & PIP (0.6) & 17.394 & 136.028 & 16.981 & 1331.195 & 14.992 & 1409.802  \\[0.02in]  \hline
         & full MIP &11.298 & 19.000 & 16.146 & 574.000 & 13.734 & 54.000  \\[0.02in]  
        Dataset 5 & full MIP &13.083 & 35.000 & 17.797 & 1097.000 & 13.734 & 3600.000  \\[0.02in]  
        ~ & PIP (0.4) & 13.080 & 19.499 & 18.103 & 343.598 & 14.859 & 1619.978  \\[0.02in]  
        ~ & PIP (0.6) & 13.085 & 31.185 & 18.103 & 1011.516 & 14.859 & 1649.590  \\[0.02in]  \hline
    \end{tabular}
\end{table}

\section{Acknowledgements}
The research of the second author is supported by the National 
Natural Science Foundation of China under grant 20221380005. 
The work of the third author was based on research supported 
by the U.S.\ Air Force Office of Scientific Research under 
grant FA9550-22-1-0045.

\newpage 
\bibliographystyle{informs2014} 

\bibliography{references}

\gap

\noindent {\bf Appendix~A:}  Consider a generic constraint 
$\phi(x) \geq 0$ where $\phi$ is a multivariate piecewise affine 
function.  By
\cite{Scholtes02}, we can write $\phi$ as the difference of two 
convex piecewise affine functions each of which is in turn the 
pointwise maximum of finitely many affine functions:
\[
\phi(x) \, = \, \displaystyle{
\max_{1 \leq i \leq I}
} \, [ \, ( a^i )^{\top}x + \alpha_i \, ] - \displaystyle{
\max_{1 \leq j \leq J}
} \, [ \, ( b^{\, j} )^{\top}x + \beta_j \, ], 
\epc x \, \in \, \mathbb{R}^n.
\]
Letting $t \triangleq \displaystyle{
\max_{1 \leq i \leq I}
} \, [ \, ( a^i )^{\top}x + \alpha_i \, ]$ be an auxiliary 
variable and
restricting $x$ to be in a bounded region so that for some constant
$\overline{C}$, we have 
$t - [ \, ( a^i )^{\top}x + \alpha_i \, ] \leq 
\overline{C}$ for all $x$ of interest, the constraint $\phi(x) \geq 0$
is then equivalent to the convex constraint 
\[
t \, \geq \,  \underbrace{\max\left\{ \, \displaystyle{
\max_{1 \leq i \leq I}
} \, [ \, ( a^{\, i} )^{\top}x + \alpha_i \, ], \ 
\displaystyle{
\max_{1 \leq j \leq J}
} \, [ \, ( b^{\, j} )^{\top}x + \beta_j \, ] \, 
\right\}}_{\mbox{convex in $x$}}
\]
subject to 
\[ \begin{array}{l}
t - ( a^{\, i} )^{\top}x - \alpha_i \, \leq \, \overline{C} v_i,
\epc i \, = \, 1, \cdots, I \\ [0.1in]
\displaystyle{
\sum_{i=1}^I
} \, v_i \, \leq \, I - 1, \\ [0.2in]
\mbox{and} \epc v_i \, \in \, \{ \, 0,1 \, \}, \ 
i \, = \, 1, \cdots, I.
\end{array}
\]
Employing this IP formulation for each constraint 
$\phi_{ij}(x) \geq \underline{B} ( 1 - z_{ij} )$ and assuming that
$\theta$ is a concave function, we obtain the following 
concave (maximizing) integer programming formulation 
of (\ref{eq:all IP}):
\[
\begin{array}{ll}
\displaystyle{
\operatornamewithlimits{\mbox{\bf maximize}}_{
x \, \in \, P; \, z; \, t; \, v}
} & \theta(x) + \displaystyle{
\sum_{j=1}^{J_0}
} \, \psi_{0j} \, z_{0j} \\ [0.2in]
\mbox{\bf subject to} & \displaystyle{
\sum_{j=1}^{J_i}
} \, \psi_{ij} \, z_{ij} \, \geq \, b_i,
\epc \forall \ i \, = \, 1, \cdots, m \\ [0.25in]
\mbox{\bf and} & \mbox{for all } j \, = \, 1, \cdots, J_i, \ 
i \, = \, 0, 1, \cdots, I \\ [0.1in]
& t_{ij} - \max\left\{ \, \displaystyle{
\max_{1 \leq \ell \leq L_{ij}}
} \, [ \, ( a_{\, ij}^{\ell} )^{\top}x + \alpha^{\ell}_{ij} \, ], \ 
\displaystyle{
\max_{1 \leq k \leq K_{ij}}
} \, [ \, ( b_{\, ij}^k )^{\top}x + \beta^k_{ij} \, ] 
\, \right\} \, \geq \, \underline{B} \, ( \, 1 - z_{ij} \, ), 
\\ [0.25in]
& t_{ij} - ( a_{\, ij}^{\ell} )^{\top}x - \alpha^{\ell}_{ij} \, \leq \, 
\overline{C} \, v_{ij}^{\ell}, \ v^{\ell}_{ij} \, \in \, \{ 0,1 \}, 
\epc \ell \, = \, 1, \cdots, L_i \\ [0.1in]
& \displaystyle{
\sum_{\ell=1}^{L_{ij}}
} \, v^{\ell}_{ij} \, \leq \, L_{ij} - 1 \ \mbox{\bf and } \ 
z_{ij} \, \in \, \{ \, 0,1 \, \},
\end{array} 
\]
where there are two sets of integer variables 
$\left\{ \, \left\{ \, \{ v^{\ell} \}_{L_{ij}}, \, z_{ij} \, 
\right\}_{j=1}^{I_i} \, \right\}_{i=0}^m$ to address the 
difference-of-convexity of the functions $\phi_{ij}$ and the 
Heaviside functions, respectively.

\gap

As an alternative to the integer formulation of the piecewise affine
constraint $\phi_{ij}(x) \geq 0$, one could treat such a constraint
directly as a difference-of-convex constraint and attempt to
apply the surrogation method in
\cite{PangRazaAlvarado16} that needs to be extended to handle
the integer variables $z_{ij}$.  Details of this extension are
outside the scope of the present paper.

\end{document}